\documentclass[11pt]{article}

\usepackage{color}
\usepackage{latexsym,dsfont}
\usepackage{amssymb, stmaryrd}
\usepackage{amsmath}
\usepackage{amsthm}
\usepackage{enumerate}
\usepackage{hyperref}
\usepackage{graphicx,float}
\usepackage{caption,mathtools}

\newtheorem{Theorem}{Theorem}[part]
\newtheorem{Definition}{Definition}[part]
\newtheorem{Proposition}{Proposition}[part]

\newtheorem{Lemma}{Lemma}[part]

\def \N{\mathbb{N}}
\def \R{\mathbb{R}}

\def \E{\mathbb{E}}
\def \F{\mathbb{F}}

\def \P{\mathbb{P}}

\def \Ac{{\cal A}}

\def \Bc{{\cal B}}

\def \Dc{{\cal D}}

\def \Fc{{\cal F}}

\def \Hc{{\cal H}}
\def \Kc{{\cal K}}
\def \Lc{{\cal L}}

\def \Oc{{\cal O}}
 \def \Nc{{\cal N}}

\def \Tc{{\cal T}}

\def \Zc{{\cal Z}}

\def \1{{\mathds 1}}

\def \eps{\varepsilon}

\def \ni{\noindent}
\def \ep{\hbox{ }\hfill$\Box$}
\def\reff#1{{\rm(\ref{#1})}}

\newcommand{\nc}{\newcommand}
\nc{\esssup}{\mathop{\mathrm{ess\,sup}}}
\nc{\essinf}{\mathop{\mathrm{ess\,inf}}}

\def\beqs{\begin{eqnarray*}}
\def\enqs{\end{eqnarray*}}
\def\beq{\begin{eqnarray}}
\def\enq{\end{eqnarray}}

\begin{document}

\title{Optimal Exploitation of a Resource \\ with Stochastic Population Dynamics\\ and Delayed Renewal}

\author{Idris Kharroubi \footnote{The research of the author benefited from the support of the French ANR research grant LIQUIRISK} \\
\footnotesize{Sorbonne Universit\'e,} \\
\footnotesize{LPSM} \\
\footnotesize{CNRS,  UMR 8001}\\
\footnotesize{\texttt{idris.kharroubi @ upmc.fr}}\and
Thomas Lim \footnote{The research of the author benefited from the support of the ``Chaire March\'e en mutation'', F\'ed\'eration Bancaire Fran\c caise}\\\footnotesize{ENSIIE,}\\
\footnotesize{LaMME} \\
\footnotesize{CNRS  UMR 8071} \\
\footnotesize{\texttt{lim @ ensiie.fr} }
\and
Vathana Ly Vath \footnote{The research of the author benefited from the support of the ``Chaire March\'e en mutation'', F\'ed\'eration Bancaire Fran\c caise}\\\footnotesize{ENSIIE,}\\
\footnotesize{LaMME} \\
\footnotesize{CNRS  UMR 8071} \\
\footnotesize{\texttt{lyvath @ ensiie.fr} }}

\maketitle

\begin{abstract}
In this work, we study the optimization problem of a renewable
resource in finite time. The resource is assumed to evolve
according to a logistic stochastic differential equation. The
manager may harvest partially the resource at any time and sell it
at a stochastic market price. She may equally decide to renew part
of the resource but uniquely at deterministic times. However, we
realistically assume that there is a delay in the renewing order.
By using the dynamic programming theory, we may obtain the PDE
characterization of our value function. To complete our study, we
give an algorithm to compute the value function and optimal
strategy. Some numerical illustrations will be equally provided.

\end{abstract}

\vspace{7mm}


\vspace{5mm}

\noindent {\bf MSC Classification (2010)~: 93E20, 62L15, 49L20, 49L25, 92D25}

\newpage


\section{Introduction}

The management of renewable resources is fundamental for the
survival and growth of the human population. An excessive
exploitation of such resources may lead to their extinction and
may therefore affect the economies of depending populations with,
for instance, high increases of prices and higher uncertainty on
the future. The typical examples are fishery \cite{C90, GKL05,
LW76} or forest management  \cite{AOK09, CR89}. Most early studies
in fishery or forest management were mainly focusing on
identifying the optimal harvesting policy. In forest economics
literature, it may be illustrated by the well-known
``tree-cutting'' problem. The most basic ``tree-cutting'' problem
is about identifying the optimal time to harvest a given forest.
Studies extending this initial tree-cutting problem have been
carried by many authors. We may, for instance, refer to
\cite{CR89} and \cite{RC90}, where the authors investigate both
single and ongoing rotation problems under stochastic prices and
forest's age or size. Rotation problem means once all the trees
are harvested, plantation takes place and planted trees may grow
up to the next harvest. In terms of mathematical formulation,
rotation problem may be reduced to an iterative optimal stopping
problem. In \cite{LSHA07}, the authors go a step further by
studying optimal replanting strategy. To be more precise, they
analyze optimal tree replanting on an area of recently harvested
forest land. However, the attempt to incorporate replanting policy
in the study of tree-cutting problem remains relatively very few,
especially when delay has to be taken into account. Indeed, the
renewed resources need some delay to become available for
harvesting. There is also an uncertainty on the renewed
quantities. In other words, the resource obtained after a renewing
decision may differ from the expected one due to some losses.
To our knowledge, these above aspects are not taken into account
in the existing literature on renewable resources management. The
aim of this paper is precisely to provide a more realistic model
in the study of optimal exploitation problems of renewable
resources by taking into account all the above features.

We suppose that the resource population evolves according to a
stochastic logistic diffusion model. Such a logistic dynamics is
classic in the modelling of populations evolution. The stochastic
aspect allows us to take into account the uncertainties of the
evolution. Since the interventions of the manager are not
continuous in practice, we consider a stochastic impulse control
problem on the resource population. We suppose that the operator
has the ability to act on the resource population through two
types of interventions. First, the manager may decide to harvest
the resource and sell the harvested resource at a given exogenous
market price. The second kind of intervention consists in renewing
the resource. Due to physical or biological constraints, the
effect of renewing orders may have some delay, i.e. a lag between
the times at which renewing decisions are taken and the time at
which renewed quantities appear in the global inventory of the
available resources. Renewing or harvesting orders are assumed
to carry both fixed and proportional costs.

From a mathematical point of view, control problems with delay
have been studied in \cite{BP09} and  \cite{OS08}, where all
interventions are delayed. Our model may be considered as more
general since some interventions are delayed while some others are
not. Another novelty of our model is the state constraints.
Indeed, the level of owned resource is a physical quantity, and
hence cannot be negative. Control problems under state
constraints, but without delay, have been studied in the
literature, see for instance \cite{MLP07} for the study of optimal
portfolio management under liquidity constraints. To deal with
such problems, the usual approach is to consider the notion of
constrained viscosity solutions introduced by Soner in
\cite{S86I,S86II}. This definition means that the value function
associated to the constrained problem is a viscosity solution in
the interior of the domain and only a semi-solution on the
boundary. In particular, the uniqueness of the viscosity solution
is usually obtained only on the interior of the domain.

In our case, we are able to characterize the behavior of the value
function on the boundary by deriving the PDE satisfied on the
frontier of the constrained domain. We therefore get the
uniqueness property of the value function on the whole closure of
the constrained domain. As a by product, we obtain the continuity
of the value function on the closure of the domain (except at
renewing dates), which improves the existing literature where this
property is obtained only on the interior of the domain, see for
instance \cite{MLP07}.

To complete our study, we provide an algorithm to compute the
value function and an associated strategy that is expected to be
optimal and apply this algorithm on a specific example.

The rest of the paper is organized as follows.  In Section
\ref{sec2}, we describe the model and the associated impulse
control problem. In Section \ref{sec3}, we give a characterization
of the value function as the unique viscosity solution to a PDE in
the class of functions satisfying a given growth condition. In
Section \ref{sec4}, we provide an algorithm to compute the value
function and an optimal strategy. Finally Section \ref{sec5} is
devoted to the proof of the main results.

%
%
%
%
%
%
%

\section{Problem formulation}\label{sec2}

\subsection{The control problem}
Let $(\Omega, \Fc, \P)$ be a complete probability space, equipped
with two mutually independent one-dimensional standard Brownian
motions $B$ and $W$. We denote by $\F:= (\Fc_t)_{t \geq 0}$ the
right-continuous and complete filtration generated by $B$ and $W$.

We consider a manager who owns a field of some given resource,
which may be exploited up to a finite horizon time $T>0$. The aim
of the manager is to manage optimally this resource in order to
maximize the expected terminal wealth which may be extracted.

\vspace{2mm}

In resource management, the manager may decide to either harvest
part of the resource or renew it. Resource renewal may be done
only at discrete times $(t_i)_{1 \leq i \leq n}$ with $t_i=i{T\over
n}$, where $n \in \N^*$. We consider an impulse control strategy $\alpha = (t_i,
\xi_i)_{1 \leq i \leq n} \cup (\tau_k, \zeta_k)_{ k \geq 1}$ where
 \begin{itemize}

\item  $\xi_i$,  $1 \leq i \leq n$, is an $\Fc_{t_i}$-measurable
random variable valued in a compact set $[0,K]$, with $K$ being a
positive constant, and corresponds to the  maximal quantity of
resource that the manager can renew,

\item $(\tau_k)_{k\geq 1}$ a nondecreasing finite or infinite
sequence of $\F$-stopping times representing the harvest times
before $T$,

\item $\zeta_k$, $k \geq 1$, an $\Fc_{\tau_k}$-measurable random
variable, valued in $\R_+$, corresponding to the harvested
quantity of resource at time $\tau_k$.
\end{itemize}

We assume the quantity of resource renewed at time $t_i$ cannot be
harvested before time $t_i + \delta$ for any $1 \leq i \leq n$
where $\delta=m{T\over n}$ with $m$ a nonnegative integer. We
suppose that for a given quantity $\xi_i$ of resource renewed at
time $t_i$, the manager may get an additional $g(\xi_i)$
harvestable resource at time $t_i+\delta=t_{i+m}$, with $g$ being
a function satisfying the following assumption.

\ni  \textbf{(H$g$)}   $g:~\R_+\rightarrow\R_+$ is a nondecreasing
and Lipschitz continuous function: there exists a positive
constant $L$ such that \beqs
|g(x)-g(x')| & \leq & L|x-x'| \;,
\enqs
for all $x,x'\in\R_+$.

\vspace{2mm}

For a given strategy $\alpha=(t_i, \xi_i)_{1 \leq i \leq n} \cup
(\tau_k, \zeta_k)_{ k \geq 1}$, we denote by $R^\alpha_t$ the
associated size of resource which is available for harvesting at
time $t$. When no intervention of the manager occurs, the
evolution of the process $R^\alpha$ is assumed to follow the below
 logistic stochastic differential equation 
\beq\label{edsR}
dR^\alpha_t & = & \eta R^\alpha_t(\lambda-
R^\alpha_t)dt+\gamma R^\alpha_t dB_t\;,
\enq
where $\eta$,
$\lambda$ and $\gamma$ are three positive constants. Since at each
time $\tau_k$, the quantity $\zeta_k$ is harvested we have
\beqs
R^\alpha_{\tau_k} & = & R^\alpha_{\tau_k^-}-\zeta_k\;.
\enqs

Moreover, we suppose that there is a natural renewal of the
resource at each time $t_i$ of a deterministic quantity
$g_0\geq0$.  Since the renewed quantity $\xi_i$ at time $t_i$ only
appears in the total resource at time $t_i+\delta=t_{i+m}$ and
increases this one of $g(\xi_i)$, we have
\beqs
 R^\alpha_{t_{i}} & = & R^\alpha_{t_{i}^-}+g_0+g(\xi_{i-m}) \;,
 \enqs
for $i=m + 1,\ldots,n$, and
 \beqs
 R^\alpha_{t_i} & = & R^\alpha_{t_i^-}+g_0 \;,
 \enqs
for $i=1,\ldots,m$.
 \vspace{2mm}

\ni The process $R^\alpha$ is then given by
\beq \nonumber
R^\alpha_t  & = &  R_0+\int_0^t\eta R^\alpha_s(\lambda- R^\alpha_s)ds+\int_0^t\gamma R^\alpha_s dB_s\\
& &  - \sum_{k \geq 1} \zeta_k \1_{\tau_k \leq t} + \sum_{i = 1}^n
g(\xi_i) \1_{t_{i + m} \leq t}+ g_0 \sum_{i = 1}^n \1_{t_i \leq
t}\;,\quad t\geq 0\;. \label{eq croissance 1}
 \enq


\vspace{2mm} We assume that the price $P$ by unit of the resource
is governed by the following stochastic differential equation \beq
\label{eq prix 1} P_t &=&P_0+\int_0^t\mu P_u du + \int_0^t\sigma
P_u dW_u \;, \quad t \geq 0 \;, \enq with $\mu$ and $\sigma$ two
positive constants.

We also define $Q_t$ the cost at time $t$ to renew a unit of the
resource. We suppose that it follows the below stochastic
differential equation \beq \label{eq prix 2} Q_t &=& Q_0 +
\int_0^t\rho Q_udu + \int_0^t\varsigma Q_udW_u \;, \quad t \geq 0
\;, \enq
  where $\rho$ and $\varsigma$ are two positive constants.

\vspace{2mm}

For a given strategy $\alpha=(t_i, \xi_i)_{1 \leq i \leq n} \cup
(\tau_k, \zeta_k)_{ k \geq 1}$, there are several costs that the
manager has to face.
\begin{itemize}
\item At each  time $\tau_k$, the manager has to pay a  cost
$c_1\zeta_k+c_2$ to harvest the quantity $\zeta_k$, where $c_1$
and $c_2$ are two positive constants. As such, by selling the
harvested quantity $\zeta_k$ at price $P_{\tau_k}$, she may get
$(P_{\tau_k} - c_1) \zeta_k - c_2$ at time $\tau_k$.

\item To renew quantity $\xi_i$ of resource at time $t_i$, the
manager has to pay  $(Q_{t_i} + c_3)\xi_i$, where $c_3$ is a
positive constant.
\end{itemize}

\vspace{2mm}

Given a control $\alpha=(t_i, \xi_i)_{1 \leq i \leq n} \cup
(\tau_k, \zeta_k)_{ k \geq 1}$ and an initial wealth $X_0$, the
wealth process $X^\alpha $ may be expressed as follows
 \beqs X_t^\alpha &=& X_0+\sum_{k \geq 1} \big[ (P_{\tau_k} -
c_1) \zeta_k - c_2 \big] \1_{\tau_k \leq t}
 - \sum_{i = 1}^n (Q_{t_i} + c_3) \xi_i \1_{t_i \leq t}\;. 
\enqs

We define the set $\Ac$ of admissible controls as the set of
strategies $\alpha$ such that \beq\label{admissiblity constraints}
\E\Big[ (X_T^\alpha)^- \Big] ~<~+\infty  &\text{ and }& R^\alpha_t
\geq 0 \quad  \text{for } 0 \leq t \leq T\;, \enq where $(.)^-$
denotes the negative part. We note that for $R_0\geq0$, the set
$\Ac$ is nonempty as it contains the strategy with no
intervention.

\vspace{2mm} We denote by $\Zc$ the set
$\Zc:=\R\times\R_+\times\R_+^*\times\R_+^*$. We define the
liquidation function $L:~\Zc \rightarrow\R$ by \beqs L(z) &: = &
\max\{x+(p-c_1)r-c_2,x\} \;,\quad \mbox{ for }~z=(x,r,p,q)\in
\Zc\;. \enqs

\ni From condition \reff{admissiblity constraints}, the
expectation $\E[L(X_T^\alpha,R_T^\alpha,P_T,Q_T)]$ is well defined
for any $\alpha\in \Ac$. We can therefore consider the objective
of the manager which consists in computing the optimal value
\beq\label{pb initial}
V_0 &:=& \sup_{\alpha \in \Ac} \E \big[L(X_T^\alpha,R_T^\alpha,P_T,Q_T) \big]\;,
\enq
and  finding a strategy $\alpha^*\in\Ac$ such that
\beq\label{opt strat}
V_0 &=&  \E\big[ L(X_T^{\alpha^*},R_T^{\alpha^*},P_T,Q_T) \big]\;.
\enq

\subsection{Value functions with pending orders}
In order to provide an analytic characterization of the value
function $V$ defined by the control problem \reff{pb initial}, we
need to extend the definition of this control problem to general
initial conditions. Moreover, since the renewing decisions are
delayed, we have to take into account the possible pending orders.

Given an impulse control $\alpha \in \Ac$, we notice that the
state of the system $R^\alpha$ is not only defined by its current
state value at time $t$ but also by the quantity at time $t$ of
the resource that has been renewed between $t-\delta$ and $t$.
 We therefore introduce the following definitions and notations. For any $t\in [0,T]$, we denote by $N(t)$ the number of possible renewing dates before $t$
 \beqs
 N(t) & := & \# \Big\{ i\in\{1,\ldots,n\} ~:~ t_i\leq t \Big\}\;,
 \enqs
  and by $D_t$ the set of renewing resource times and the associated quantities  between $t-\delta$ and $t$
\beq\label{def D_t} D_t &:=& \Big\{d=(t_i,e_i)_{N(t-\delta)+1 \leq
i \leq N(t)}~:~ e_i\in  \R_+ \mbox{ for }
i=N(t-\delta)+1,\ldots,N(t) \Big\} \;,\quad \enq with the
convention that $D_t= \emptyset$ if $N(t-\delta)=N(t)$.

\vspace{2mm}

For any $t \in [0,T]$ and  $d=(t_i,e_i)_{N(t-\delta)+1 \leq i \leq
N(t)} \in D_t$, we denote by $\tilde \Ac_{t,d}$ the set of
strategies which take into account the pending renewing decisions
taken between $t - \delta$ and $t$ \beqs
\tilde \Ac_{t,d} &:=& \Big\{\alpha=(t_i, \xi_i)_{N(t-\delta)+1 \leq i \leq n} \cup (\tau_k, \zeta_k)_{ k \geq 1}  ~: ~ \\
 & & \quad \xi_i=e_i ~\mbox{ for }~i=N(t-\delta)+1, \ldots , N(t)\;; \\
 & & \quad\xi_i \mbox{ is } \Fc_{t_i}-\mbox{measurable  for }~N(t)+1 \leq i \leq n \;; \\
 & & \quad (\tau_k)_{k \geq 1} \mbox{ is a nondecreasing finite or infinite sequence of } \F-\mbox{stopping time}  \mbox{ with }\tau_1 > t\;; \\
& & \quad \zeta_k \mbox{ is } \Fc_{\tau_k}-\mbox{measurable  for }~k\geq 1 \Big\}\;.
\enqs

For  $z=(x,r,p,q) \in \Zc$,   $d \in D_t$ and $\alpha \in \tilde
\Ac_{t,d}$, we denote by $Z^{t,z,\alpha}=
(X^{t,z,\alpha},R^{t,r,\alpha},P^{t,p},Q^{t,q})$ the quadruple of
processes defined by
 \beq R^{t,r,\alpha}_s &=& r+ \int_t^s \eta
R^{t,r,\alpha}_u(\lambda-R^{t,r,\alpha}_u)du + \int_t^s \gamma
R^{t,r,\alpha}_u  dB_u - \sum_{k\geq 1} \zeta_k \1_{\tau_k \leq s}
\nonumber \\\label{eq croissance 1 bis}
 & & + \sum_{i = N(t-\delta)+1}^n g(\xi_i) \1_{t_{i + m} \leq s} +g_0\big(N(s)-N(t)\big)
 \;, \\\label{dyn X}
 X^{t,z,\alpha}_s & = & x+ \sum_{k \geq 1} \big[ (P^{t,p}_{\tau_k} - c_1) \zeta_k- c_2 \big] \1_{\tau_k \leq s} - \hspace{-4mm}\sum_{i = N(t)+1}^n \hspace{-4mm}(Q^{t,q}_{t_i} +c_3 ) \xi_i  \1_{t_i \leq s}\;,\\ \label{dyn P}
  P^{t,p}_s & = & p+\int_t^s\mu P^{t,p}_u du + \int_t^s\sigma  P^{t,p}_u dW_u\;,\\\label{dyn Q}
  Q^{t,q}_s & = & q+\int_t^s\rho Q^{t,q}_u du + \int_t^s\varsigma  Q^{t,q}_u dW_u \;,
 \enq
 for $s\in[ t,T]$.
 We denote by $\Ac_{t,z,d}$ the set of strategies $\alpha\in \tilde \Ac_{t,d} $ such that
 \beq
  \E\Big[ (X^{t,z,\alpha}_{T})^-\Big]~<~+\infty &\text{and} &  R^{t,r,\alpha}_s \geq 0 
  \quad
  \text{for all } s \in[t, T]\;. \qquad \label{cond-adm}
 \enq
We then consider for $(t,z) \in [0,T] \times \Zc$,   $d \in D_t$, $\alpha \in \Ac_{t,z,d}$
the following benefit criterion
\beqs
J(t,z,\alpha) &:=& \E \Big[ L(Z^{t,z,\alpha}_T)
 \Big]\;,
\enqs
which is well defined under conditions \reff{cond-adm}.
 We define the corresponding value function by
\beqs
v(t,z,d) &:=& \sup_{\alpha \in \Ac_{t,z,d}} J(t,z,\alpha) \;, \quad  (t,z,d) \in \mathcal D\;,
\enqs
where $\mathcal D$ is the definition domain of $v$ defined by
\beqs
\Dc & = & \Big\{(t,z,d)~:~(t,z)\in[0,T]\times \Zc \mbox{ and } d\in D_t\Big\}\;.
\enqs

\vspace{2mm}

\ni For simplicity, we also introduce the operators $\Gamma^{c}$,
$\Gamma^{rn}_1$ and $\Gamma^{rn}_2$ given by \beqs
\Gamma^{c}(z,\ell) & := & (x+(p-c_1)\ell-c_2,r-\ell,p,q) \;,\\
\Gamma^{rn}_1(z,\ell) & := & (x-(q+c_3)\ell,r+g_0,p,q) \;,\\
\Gamma^{rn}_2(z,\ell) & := & (x,r+g(\ell),p,q)\;, \enqs for all
$z=(x,r,p,q)\in\Zc$ and $\ell\in\R_+$. The operator $\Gamma^c$
corresponds to the new position of the state process after a
resource consumption decision: if the manager harvests $\zeta_k$
at time $\tau_k$, then the state process is \beqs
Z^{t,z,\alpha}_{\tau_k} & = &
\Gamma^c(Z^{t,z,\alpha}_{\tau_k^-},\zeta_k) \;, \enqs and
$\Gamma^{rn}_1$ and $\Gamma^{rn}_2$ correspond to the new position
of the state process after a renewal decision: if the manager
renews $(\xi_i)_{1 \leq i \leq n}$ at times $(t_i)_{1 \leq i \leq
n}$, then the state process is given by \beqs
Z^{t,z,\alpha}_{t_i} & = & \Gamma^{rn}_1(Z^{t,z,\alpha}_{t_i^-},\xi_i)\;, \quad \mbox{for } i = 0 , \ldots, m \;,\\
Z^{t,z,\alpha}_{t_{i}} & = &
\Gamma^{rn}_1(\Gamma^{rn}_2(Z^{t,z,\alpha}_{t_{i}^-},\xi_{i-m}),\xi_i
) \;, \quad \mbox{for } i = m+1 , \ldots, n \;.
\enqs

We first give a new expression of the value function $v$. To this
end, we introduce the set \beqs \hat \Ac_{t,z,d}  & = &
\Big\{\alpha=(t_i, \xi_i)_{N(t-\delta)+1 \leq i \leq n} \cup
(\tau_k, \zeta_k)_{ k \geq 1}\in \tilde \Ac_{t,d}  ~: \\
 & &
~\big(P^{t,p}_{\tau_k}-c_1\big)\zeta_k-c_2~\geq~ 0~ \forall
\,k\geq 1~ \mbox{Êand } R^{t,r,\alpha}_s \geq 0~ \forall
\,s\in[t,T] \Big\}\;. \enqs

\begin{Proposition}
The value function $v$ can be expressed as follows
\beq\label{eq vA hat}
v(t,z,d) & = & \sup_{\alpha \in \hat \Ac_{t,z,d}} J(t,z,\alpha) \;, \quad  (t,z,d) \in \mathcal D\;.
\enq
\end{Proposition}
\ni \textbf{Proof.} Fix $(t,z,d) \in \mathcal D$ with $z=(x,r,p,q)$ and denote by $\hat v(t,z,d) $ the right hand side of \reff{eq vA hat}.

We first notice that  $\hat \Ac_{t,z,d}\subset  \Ac_{t,z,d}$. Indeed, for $\alpha \in \hat \Ac_{t,z,d}$, we have
\beqs
X^{t,z,\alpha}_T & = & x+\sum_{k \geq 1} \big[ (P^{t,p}_{\tau_k} - c_1) \zeta_k- c_2 \big] \1_{\tau_k \leq T} - \hspace{-4mm}\sum_{i = N(t)+1}^n \hspace{-4mm}(Q^{t,q}_{t_i} +c_3 ) \xi_i  \\
 & \geq & x-nK\big(\sup_{s\in[t,T]}Q^{t,q}_s+ c_3\big)\;.
\enqs
Since $Q^{t,q}$ follows the dynamics \reff{dyn Q}, we have $\E[\sup_{s\in[t,T]}Q^{t,q}_s]<+\infty$ and we get $\E[(X^{t,z,\alpha}_T)^-]  <  +\infty$. We therefore deduce that
\beqs
v(t,z,d) & \geq & \hat v(t,z,d)\;.
\enqs
We turn to the reverse inequality. Fix $\alpha=(t_i, \xi_i)_{N(t-\delta)+1 \leq i \leq n} \cup (\tau_k, \zeta_k)_{ k \geq 1}\in \Ac_{t,z,d}$ and define the associated strategy $\hat \alpha=(t_i, \xi_i)_{N(t-\delta)+1 \leq i \leq n} \cup (\hat \tau_k, \hat \zeta_k)_{ k \geq 1}\in \hat \Ac_{t,z,d}$ by
\beqs
(\hat \tau_j,\hat \zeta_j) & = & (\tau_{k_j},\zeta_{k_j})\quad ~~\mbox{ for } j\geq 1\;,
\enqs
where the sequence $(k_j)_{j\geq 1}$ is defined by
\beqs
k_1 & = & \min \{k\geq 1~:~(P^{t,p}_{\tau_k}-c_1)\zeta_k-c_2\geq 0\}\;,\\
k_j & = & \min\{k\geq k_{j-1}+1~:~(P^{t,p}_{\tau_k}-c_1)\zeta_k-c_2\geq 0\}\;,
\enqs
$i.e.$ $\hat \alpha$ is obtained from $\alpha$ by keeping only harvesting orders such that $(P^{t,p}_{\tau_k}-c_1)\zeta_k-c_2\geq 0$.

We then easily check from dynamics \reff{eq croissance 1 bis} and \reff{dyn X} that
\beqs
X_s^{t,z,\alpha} ~ \leq ~ X_s^{t,z,\hat \alpha} & \mbox{ and } & R_s^{t,r,\alpha} ~ \leq ~ R_s^{t,r,\hat \alpha}
\enqs
for all $s\in[t,T]$.
Therefore we get
\beqs
L(Z_T^{t,z,\alpha}) & \leq & L(Z_T^{t,z,\hat \alpha})\;,
\enqs
which gives
\beqs
 \hat v(t,z,d) & \geq & v(t,z,d)\;.
\enqs
\ep
\section{PDE characterization}\label{sec3}

\subsection{Boundary condition and dynamic programming principle}


 We first provide a boundary condition for the value function associated to the optimal management of renewable resource.

 \vspace{2mm}

 \begin{Proposition}\label{Bound value function}
 The value function $v$ satisfies the following growth condition: there exists a constant $C$ such that
 \beq\label{growth-property}
x
~~\leq~~ v(t,z,d) & \leq & x+C\Big(1+|r|^{4}+|p|^{4}+|q|^{4}\Big) \;,\qquad
 \enq
 for all $t\in[0,T]$, $z=(x,r,p,q)\in\Zc$,  and $d\in D_t$.
 \end{Proposition}
 \ni The proof of this proposition is postponed to Section \ref{preuve1}.
 \vspace{2mm}

With this bound, we are able to state the dynamic programming
relation on the value function of our control problem with
execution delay. For any $t\in [0,T]$, $d\in D_t$ and $\alpha
=(t_i,\xi_i)_{N(t-\delta)+1\leq i\leq n}\cup(\tau_k,\zeta_k)_{k
\geq 1} \in \hat \Ac_{t,z,d}$, we denote \beqs
d(u,\alpha)&=& (t_{i}, \xi_{i})_{N(u-\delta)+1 \leq i \leq N(u)}\;,\quad u\in[t,T] \;,
\enqs
with the convention that $d(u,\alpha)= \emptyset$ if $N(u-\delta)=N(u)$.
We notice that $d(u,\alpha)$ corresponds to the set of renewing orders that have been given before $u$ and whose delayed effects appear after $u$. We also denote by $\Tc_{[t,T]}$ the set of $\F$-stopping times valued in $[t,T]$.
\begin{Theorem}\label{DP}
The value function $v$ satisfies the following dynamic programming principle.
%
%


 \begin{enumerate}[(DP1)]
\item First dynamic programming inequality:
\beqs
  v(t,z,d)& \geq & \E \Big[  v_{}(\vartheta,Z^{t,z,\alpha}_\vartheta,d(\vartheta,\alpha))\Big]\;,
\enqs
for all $\alpha \in \hat \Ac_{t,z,d}$ and all $\vartheta\in\Tc_{[t,T]}$.

\item Second dynamic programming inequality: for any $\varepsilon >0$, there exists $\alpha \in \hat \Ac_{t,z,d}$ such that
  \beqs
  v(t,z,d) - \varepsilon & \leq & \E \Big[ v_{}(\vartheta,Z_\vartheta^{t,z,\alpha},d(\vartheta,\alpha))\Big] \;,
\enqs
for all $\vartheta\in\Tc_{[t,T]}$.
\end{enumerate}
\end{Theorem}
\vspace{2mm}

 \ni The proof of this proposition is postponed to Section \ref{preuve2}.
 \vspace{2mm}

\subsection{Viscosity properties and uniqueness}
The PDE system associated to our control problem is formally
derived from the dynamic programming relations. We first decompose
the domain $\Dc$ as follows
\beqs \Dc & = & \bigcup_{k=0}^{n}\Dc_k
\;, \enqs where \beqs \Dc_k &= & \Big\{(t,z,d)\in \Dc
~:~~t\in\big[t_k,t_{k+1}\big)\Big\}\;, \enqs for $k=0,\ldots,n-1$
and \beqs \Dc_n & = & \Big\{(t,z,d)\in \Dc~:~t=T\Big\}\;. \enqs We
also decompose the sets $\Dc_k$, $k=0,\ldots,n$, as follows \beqs
\Dc_k & = &  \Dc_k^1\cup\Dc_k^2 \;, \enqs where \beqs
\Dc_k^1 & = & \big\{(t,z,d)\in\Dc_k~:~z=(x,r,p,q) \mbox{ with } r=0\big\}\;,\\
\Dc_k^2 & = & \big\{(t,z,d)\in\Dc_k~:~z=(x,r,p,q) \mbox{ with } r>0\big\}\;.
\enqs
 We define the operators $\Hc$, $\Nc_1$, $\bar\Nc_1$, $\Nc_2$ and $\bar \Nc_2$ by
\beqs
\Hc \phi (t,z,d) & = & \sup_{0\leq a \leq r}  \phi\big(t,\Gamma^c(z,a),d \big) \;,
\enqs
for  any $(t,z,d)\in \Dc$ and any function $\phi$ defined on $\Dc$,
\beqs
\Nc_1 \phi (t_k,z,d)  & = & \sup_{
0 \leq e \leq K} \phi \Big(t_k,\Gamma^{rn}_1 \big(\Gamma^{rn}_2(z,e_{k-m}),e \big),d\cup(t_k,e)\setminus (t_{k-m}, e_{k-m}) \Big) \;,\\
\bar  \Nc_1 \phi (t_k,z,d)  & = & \sup_{
\begin{tiny}\begin{array}{c}
0 \leq e \leq K\\
0\leq a \leq r \end{array}\end{tiny}
 } \phi \Big(t_k,\Gamma^{rn}_1 \big(\Gamma^{rn}_2\big(\Gamma^c(z,a),e_{k-m}\big),e \big),d\cup(t_k,e)\setminus (t_{k-m}, e_{k-m}) \Big) \;,
\enqs
for any $(t_k,z,d)\in \Dc$ with $k=m + 1,\ldots,n$,  and any function $\phi$ defined on $\Dc$, and
\beqs
\Nc_2 \phi (t_k,z,d)  & = & \sup_{0 \leq e \leq K
} \phi \Big(t_k,\Gamma^{rn}_1 \big(z,e \big),d\cup(t_k,e) \Big) \;,\\
\bar \Nc_2 \phi (t_k,z,d)  & = & \sup_{
\begin{tiny}\begin{array}{c}
0\leq e \leq K\\
0\leq a\leq r \end{array}\end{tiny}
} \phi \Big(t_k,\Gamma^{rn}_1 \big(\Gamma^c(z,a),e \big),d\cup(t_k,e) \Big) \;,
\enqs
for any $(t_k,z,d)\in \Dc$ with $k=0,\ldots,m$, and any function $\phi$ defined on $\Dc$.

\vspace{2mm}

This provides equations for the value function $v$ which takes the
following nonstandard form \beq\label{EDP0}
 - \Lc v(t,z,d) & = & 0
\enq
for $(t,z,d)\in \Dc^1_k$, with $k=0,\ldots,n$,
\beq\label{EDP1}
\min\Big\{ - \Lc v(t,z,d) ~,~
v(t,z,d) - \Hc v(t,z,d)
\Big\} &=&0
\enq
for $(t,z,d)\in \Dc^2_k$, with $k=0,\ldots,n-1$,
\beq
v(T^-,z,d) & = & \max\big\{\Nc_1 L(z,d)~,~\bar \Nc_1 L(z,d)\big\}\qquad\label{EDP2}
\enq
for $(T,z,d)\in \Dc$,
 \beq\label{EDP4}
v(t_{k}^-,z,d) & = & \max\{ \Nc_1 v(t_k,z,d)~,~\bar \Nc_1v(t_k,z,d)\}
\enq
for  $(t_k,z,d)\in \Dc_k$, with $k=m+1,\ldots,n-1$,
and
\beq
v(t_{k}^-,z,d) & = &  \max\big\{\Nc_2 v(t_k,z,d)~,~\bar \Nc_2 v(t_k,z,d)\big\}\qquad\label{EDP6}
\enq
for $(t_k,z,d)\in \Dc_k$, with $k=0,\ldots,m$.

\vspace{2mm}

Here $\Lc$ is the second order local operator associated to the diffusion $(P,Q,R)$ with no intervention. It is given by
\beqs
\Lc \varphi (t,z) & = & \partial_t \varphi(t,z) + \mu p \partial_p\varphi(t,z)+\rho q\partial_q\varphi(t,z)+\eta r(\lambda-r)\partial_r\varphi(t,z) \\
 &  & +{1\over 2} \Big(\sigma^2p^2\partial^2_{pp}\varphi(t,z)+\varsigma^2q^2\partial^2_{qq}\varphi(t,z)+2\sigma\varsigma pq \partial ^2_{pq}\varphi(t,z)+\gamma^2r^2 \partial ^2_{rr}\varphi(t,z)\Big)
\enqs
for any $(t,z)\in[0,T]\times\Zc$ with $z=(x,r,p,q)$ and any function $\varphi\in C^{1,2}([0,T]\times \Zc )$.

\vspace{2mm}

As usual, we do not have any regularity property on the value
function $v$. We therefore work with the notion of (discontinuous)
viscosity solution. Since our system of PDEs \reff{EDP0} to
\reff{EDP6} is nonstandard, we have to adapt the definition to our
framework.

\vspace{2mm}

First, for a locally bounded function $w$ defined on $\Dc$, we define its lower semicontinuous (resp. upper semicontinuous) envelop $w_*$ (resp. $ w^*$) by
\beqs
w _* (t,z,d) & = & \liminf_{\begin{tiny}\begin{array}{c}(t',z',d')\rightarrow(t,z,d)\\(t',z',d')\in \Dc_k\end{array}\end{tiny}} w(t',z',d') \;,\\
w ^* (t,z,d) & = & \limsup_{\begin{tiny}\begin{array}{c}(t',z' ,d')\rightarrow(t,z,d)\\(t',z',d')\in \Dc_k\end{array}\end{tiny}} w(t',z',d') \;,\enqs
for $(t,z,d)\in \Dc_k $, with $k=0,\ldots,n-1$. We also define its left lower semicontinuous (resp. upper semicontinuous) envelop at time $t_k$ by
\beqs
w _* (t_k^-,z,d) & = & \liminf_{\begin{tiny}\begin{array}{c}(t',z',d')\rightarrow(t_k^-,z,d)\\(t',z',d')\in \Dc_{k-1}\end{array}\end{tiny}} w(t',z',d)\;,\\
w ^* (t_k^-,z,d) & = & \limsup_{\begin{tiny}\begin{array}{c}(t',z',d')\rightarrow(t_k^-,z,d)\\(t',z',d')\in \Dc_{k-1}\end{array}\end{tiny}} w(t',z',d) \;,
\enqs
for $k \in \{1,\ldots,n\}$.
\begin{Definition}[Viscosity solution to \reff{EDP0} -- \reff{EDP6}]
 A  locally bounded function $w$ defined on $\Dc$ is a viscosity supersolution (resp. subsolution) if

\ni (i) for any $k=0,\ldots,n-1$, $(t,z)\in\Dc_k^1$ and $\varphi\in C^{1,2}(\Dc_k)$ such that
\beqs
(w_*-\varphi)(t,z,d) & = & \min_{\Dc_k} (w_*-\varphi)\\
\mbox{(resp. } ( w ^*-\varphi)(t,z,d) & = & \max_{\Dc_k} (w ^*-\varphi)\mbox{)}
\enqs
we have
\beqs
- \Lc \varphi(t,z,d) &\geq  &0\\
\mbox{(resp. }  - \Lc \varphi(t,z,d) & \leq & 0 \mbox{)}\;,
\enqs
(ii) for any $k=0,\ldots,n-1$, $(t,z)\in\Dc^2_k$ and $\varphi\in C^{1,2}(\Dc_k)$ such that
\beqs
(w_*-\varphi)(t,z,d) & = & \min_{\Dc_k} (w_*-\varphi)\\
\mbox{(resp. } ( w ^*-\varphi)(t,z,d) & = & \max_{\Dc_k} (w ^*-\varphi)\mbox{)}
\enqs
we have
\beqs
\min\Big\{ - \Lc \varphi(t,z,d) ~,~
w_*(t,z,d) - \Hc w_*(t,z,d)
\Big\} &\geq  &0\\
\mbox{(resp. } \min\Big\{ - \Lc \varphi(t,z,d) ~,~
w^*(t,z,d) - \Hc w^*(t,z,d)
\Big\} & \leq & 0 \mbox{)}\;,
\enqs
(iii) for any $(T,z,d)\in \Dc$ we have
\beqs
w_*(T^-,z,d) & \geq   & \max\{ \Nc_1 L(z,d)~,~\bar \Nc_1 L(z,d)\}\\
\mbox{(resp. } w_*(T^-,z,d) & \leq & \max\{\Nc_1 L(z,d)~,~\bar \Nc_1 L(z,d)\} \mbox{)}\;,
\enqs
(iv) for any $k=m+1,\ldots,n-1$, $(t_k,z,d)\in \Dc$ we have
\beqs
w_*(t_{k}^-,z,d) & \geq & \max\{ \Nc_1 w_*(t_k,z,d)~,~\bar \Nc_1w_*(t_k,z,d)\}\\
\mbox{(resp. } w^*(t_{k}^-,z,d) & \leq & \max\{ \Nc_1 w^*(t_k,z,d)~,~\bar \Nc_1w^*(t_k,z,d)\}\mbox{)}\;,
\enqs
(v) for any $k=0,\ldots,m$, $(t_k,z,d)\in \Dc$  we have
\beqs
w_*(t_{k}^-,z,d) & \geq & \max\{ \Nc_2 w_*(t_k,z,d)~,~\bar \Nc_2w_*(t_k,z,d)\}\\
\mbox{(resp. } w^*(t_{k}^-,z,d) & \leq & \max\{ \Nc_2 w^*(t_k,z,d)~,~\bar \Nc_2w^*(t_k,z,d)\}\mbox{)}\;.
\enqs

\vspace{2mm}

A  locally bounded function $w$ defined  on $\Dc$ is said to be a
viscosity solution to \reff{EDP0}--\reff{EDP6} if it is a
supersolution and a subsolution to \reff{EDP0}--\reff{EDP6}.

\end{Definition}

\ni The next result provides the viscosity properties of the value function $v$.
\begin{Theorem}[Viscosity characterization]
The value function $v$ is the unique viscosity solution to
\reff{EDP0}--\reff{EDP6} satisfying the growth condition
\reff{growth-property}. Moreover, $v$ is continuous on $\Dc_k$ for
all $k=0,\ldots,n-1$.
\end{Theorem}


\section{Numerics}\label{sec4}

We describe, in this section, a backward algorithm to approximate
the value function and an optimal strategy. Some numerical
illustrations are also provided.

\subsection{Approximation of the value function $v$}

\paragraph{Initialization step.} For $(t,z,d)\in \Dc_{n-1}^1$ we have
\beqs
v(t,z,d) & = & \E\Big[  \max\{ \Nc_1 L(Z^{t,z,d}_T,d)~,~\bar \Nc_1 L(Z^{t,z,d}_T,d)\} \Big]\;.
\enqs
We can therefore approximate it by $\hat v(t,z,d)$ which is the associated Monte Carlo estimator.

On $\Dc_{n-1}^2$ the function $v$ is solution to the PDE
\reff{EDP1} with the terminal condition \reff{EDP2}. Therefore, we
can  compute an approximation $\hat v$   using an  algorithm
computing  optimal values of impulse control problem with boundary
on $\Dc_{n-1}^1$ and the terminal value given by  \reff{EDP2}
(see e.g. \cite{COS02}).

\paragraph{Step $k+1\rightarrow k$.} Once we have an approximation $\hat v(t,z,d)$ of $v(t,z,d)$ for
$(t,z,d)\in \Dc_{k+1}$ we are able to get an approximation of $v$
on $\Dc_{k}$ as follows.

\vspace{2mm}

\ni $\bullet$ Case 1: $m\leq k\leq n-1$.

For $(t,z,d)\in \Dc_{k}^1$ we have
\beqs
v(t,z,d) & = & \E\Big[  \max\{ \Nc_1 v(t_{k+1},Z^{t,z,d}_{t_{k+1}},d)~,~\bar \Nc_1v(t_{k+1},Z^{t,z,d}_{t_{k+1}},d)\}  \Big]\;.
\enqs
We can therefore approximate it by $\hat v(t,z,d)$ which is the  Monte Carlo estimator of
\beqs
\E\Big[  \max\{ \Nc_1\hat v(t_{k+1},Z^{t,z,d}_{t_{k+1}},d)~,~\bar \Nc_1 \hat v(t_{k+1},Z^{t,z,d}_{t_{k+1}},d)\}  \Big]\;.
\enqs
 On $\Dc_{k}^2$ the function $v$ is solution to the PDE
\reff{EDP1} with the terminal condition \reff{EDP4}.  Since we
already have approximations of $v$ on $\Dc^1_{k}$ and
$\Dc_{k+1}$, we can compute an approximation $\hat v$  using an
algorithm computing  optimal values of impulse control problem
with boundary  on $\Dc_{k}^1$  (see e.g. \cite{COS02}) and the
terminal value given by  \beqs \hat v(t_{k+1}^-,z,d) & = & \max\{
\Nc_1 \hat v(t_{k+1},z,d)~,~\bar \Nc_1\hat v(t_{k+1},z,d)\}\;.
\enqs

\vspace{2mm}

\ni $\bullet$ Case 2: $0\leq k\leq m - 1$.  The procedure is the same as in Case 1 but with $\Nc_2$ and $\bar \Nc_2$ instead of   $\Nc_1$ and $\bar \Nc_1$ respectively.

 \vspace{2mm}

\subsection{An optimal strategy for the approximated problem}

We turn to the computation of an optimal strategy. From the general optimal stopping theory (see \cite{EK81}), we provide the following strategy $\hat \alpha$. This strategy is constructed as usually done for optimal strategies of impulse control problem but using the approximation $\hat v$ instead of the value function $v$.
We start with an initial data $(t,z,d)$. We denote by $\hat \alpha=(t_i,\hat \xi_i)_{N(t-\delta)+1\leq i\leq n}\cup(\hat \tau_k,\hat \zeta_k)_{k \geq 1}$ the strategy constructed step by step and by $\hat Z^\kappa=(\hat X^\kappa,\hat R^\kappa,\hat P^\kappa,\hat Q^\kappa)$ the process controlled by the truncated strategy $\hat \alpha ^\kappa:= (t_i,\hat \xi_i)_{N(t-\delta)+1\leq i\leq n}\cup(\hat \tau_k,\hat \zeta_k)_{\kappa\geq k \geq 1}$. We also denote by $\hat d _s=(t_i,\hat e_i)_{N(s-\delta)+1\leq i\leq N(s)}$ the pending orders at time $s\in[t,T]$.

\paragraph{Initialization step.} We first start by computing the first  harvesting time $\hat \tau_1$ by
\beqs
\hat \tau_1 & = & \inf\Big\{s\geq t~:~\hat v(s, \hat Z^0_s, \hat d_s)~=~\Hc \hat v(s, \hat Z^0_s, \hat d_s) \Big\}
\enqs
and the associated harvested quantity $\hat \zeta _1$ by
\beqs
\hat \zeta_1 & \in & \text{arg}\max_{0\leq a\leq \hat R^0_{\tau_1}} \hat v(\hat \tau_1, \Gamma^c(\hat Z^0_{\hat \tau_1},a), \hat d_{\hat \tau_1})\;.
\enqs

 \paragraph{Step $k\rightarrow k+1$ for harvesting orders.} We then compute the $(k+1)$-th  harvesting time $\hat \tau_{k+1}$ by
\beqs
\hat \tau_{k+1} & = & \inf\Big\{s\geq \hat \tau_{k}~:~\hat v(s, \hat Z^k_s, \hat d_s)~=~\Hc \hat v(s, \hat Z^k_s, \hat d_s) \Big\}
\enqs
and the associated harvested quantity $\hat \zeta _{k+1}$ by
\beqs
\hat \zeta_{k+1} & \in & \text{arg}\max_{0\leq a\leq \hat R^k_{\tau_{k+1}}} \hat v(\hat \tau_{k+1} , \Gamma^c(\hat Z^{k}_{\hat \tau_{k+1} },a), \hat d_{\hat \tau_{k+1} })\;.
\enqs

 \paragraph{Step $i$ for renewing orders.} Denote by $\hat k_s$ the (random) number of harvesting orders on $[t,s]$.
 We then distinguish two cases.

 \vspace{2mm}

\ni  $\bullet$ Case 1: $0\leq i\leq m$.

\ni Suppose first that
\beqs
 \Nc_2 \hat v( t_i, \hat Z ^{\hat k_{t_i}} _{t_i},\hat d_{t_{i-1}}) & \geq & \bar \Nc_2 \hat v( t_i, \hat Z ^{\hat k_{t_i}} _{t_i},\hat d_{t_{i-1}})\;.
\enqs
 Then we compute the optimal renewed resource $\hat \xi_i$ at time $t_i$ by
 \beqs
 \hat \xi_i & = & \text{arg} \max_{0\leq e\leq K} \hat v\Big( t_i, \Gamma_1^{rn}\big(\hat Z ^{\hat k_{t_i}} _{t_i},e\big),\hat d_{t_{i-1}}\cup(t_i,e)\Big) \;.
 \enqs
 \ni If we now suppose that
\beqs
 \Nc_2 \hat v( t_i, \hat Z ^{\hat k_{t_i}} _{t_i},\hat d_{t_{i-1}}) & < & \bar \Nc_2 \hat v( t_i, \hat Z ^{\hat k_{t_i}} _{t_i},\hat d_{t_{i-1}})\;.
\enqs
 Then we compute the optimal renewed resource $\hat \xi_i$ at time $t_i$ by
 \beqs
 \hat \xi_i & = & \text{arg} \max_{0\leq e\leq K} \hat v\Big( t_i, \Gamma_1^{rn}\big(\Gamma_1^{c}\big(\hat Z ^{\hat k_{t_i^-}} _{t_i}, \hat \zeta_{\hat k_{t_i}}\big),e\big),\hat d_{t_{i-1}}\cup(t_i,e)\Big)
 \enqs
which is also given by the same expression as in the first inequality
 \beqs
 \hat \xi_i & = & \text{arg} \max_{0\leq e\leq K} \hat v\Big( t_i, \Gamma_1^{rn}\big(\hat Z ^{\hat k_{t_i}} _{t_i},e\big),\hat d_{t_{i-1}}\cup(t_i,e)\Big)\;.
 \enqs

 \vspace{2mm}

\ni  $\bullet$ Case 2: $m+1\leq i\leq n$.

 \ni As in the first case we do not need to distinguish the subcases $\Nc_1 \hat v\geq \bar \Nc_1 \hat v$ and  $\Nc_1 \hat v<\bar \Nc_1 \hat v$ and the optimal renewed quantity at time $t_i$ is given by
  \beqs
 \hat \xi_i & = & \text{arg} \max_{0\leq e\leq K} \hat v\Big( t_i, \Gamma_1^{rn}\big(\Gamma_2^{rn}\big(\hat Z ^{\hat k_{t_i}} _{t_i},\hat e_{i-m}\big),e\big),\hat d_{t_{i-1}}\cup(t_i,e)\setminus (t_{i-m},\hat e_{i-m})\Big)\;.
 \enqs

\subsection{Examples}
In this part we present numerical illustrations that we get by using an implicit finite difference scheme mixed with an iterative procedure which leads to the resolution of a Controlled Markov Chain by assuming that the resource is a forest . This class of problems is intensively studied by Kushner and Dupuis \cite{KushDup}. The convergence of the solution of the numerical scheme towards the solution of the HJB equation, when the time-space step goes to zero, can be shown using the standard local consistency argument i.e. the first and second moments of the approximating Markov chain converge to those of the continuous process $(R,P)$. We assume that the maximal size of the forest is 1 and we use a discretization step of $1/151$ for the size of the forest. About the discretization of the price we discretize the process $S = \log(P)$ with $P_0=1$, we consider $S_{\min} = - |\mu - \sigma^2/2|*T - 3 \sigma \sqrt{T}$ and $S_{\max} =  |\mu - \sigma^2/2|*T + 3 \sigma \sqrt{T}$, and the discretization step is $1/101$.

 We compute the optimal strategy to harvest and renew, and the value function. We assume the parameters of the logistic SDE are $\eta=1$, $\lambda=0.7$ and $\gamma=0.1$. The parameter of natural renewal is $g_0=3\%$ of the forest. The delay before to able to harvest a tree which is renewed is 1 and the function $g(x)$ is equal to $x$. The initial price is $1$. The parameters of the price $P$ are $\mu=0.07$ and $\sigma=0.1$, and the costs to harvest and renew are $c_1=0.1$, $c_2=0.01$ and $c_3=0.1$. We assume that the price $Q$ is equal to the price $P$. We can renew at times $\{1,2\}$ and the terminal time is $T=3$.
\begin{figure}[H]
\includegraphics[width=14cm]{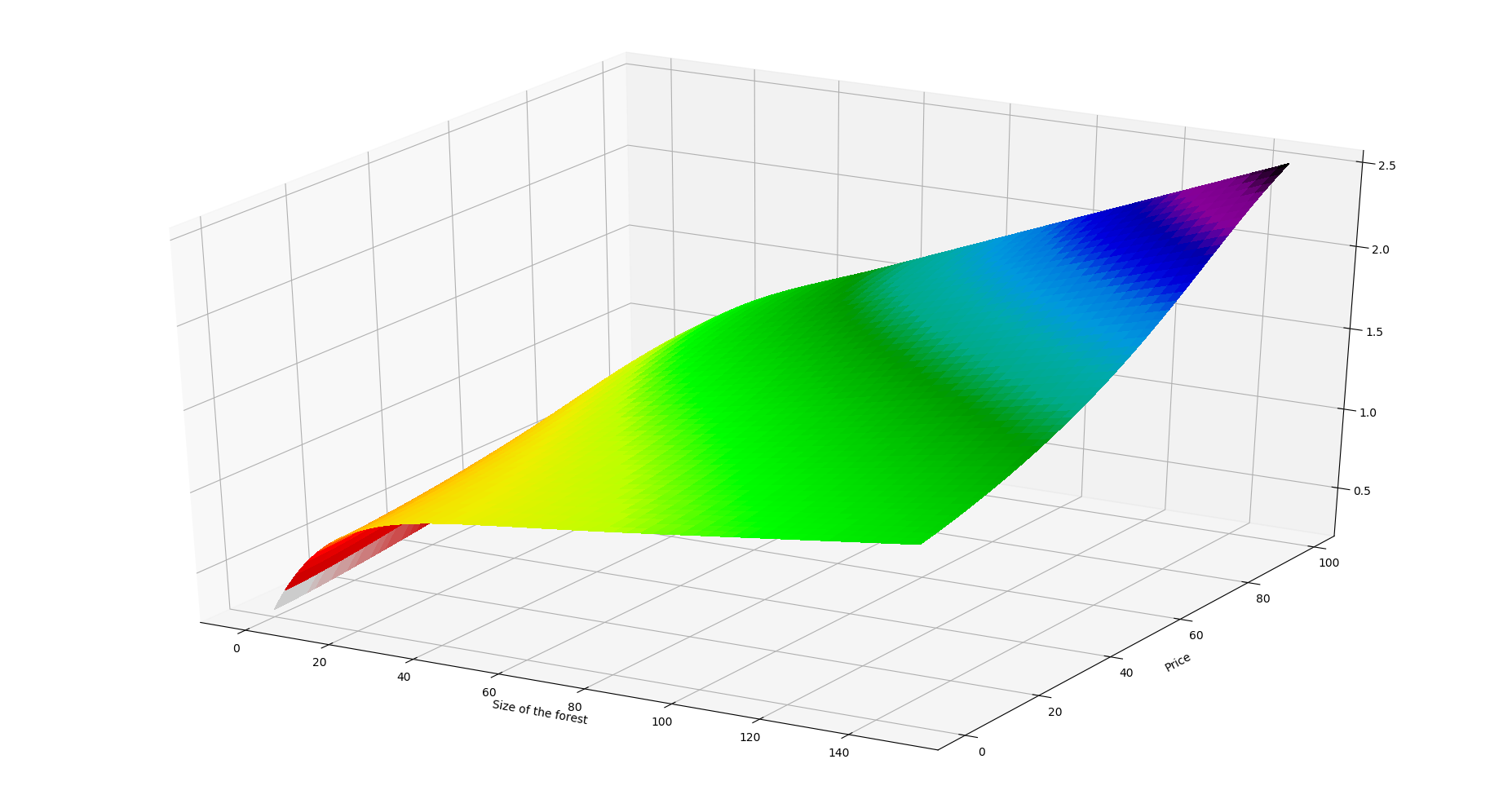}
\caption{The value function with respect to the price $P$ and the size of the forest $R$.}\label{fig:valeur}
\end{figure}
We remark that the value function is increasing w.r.t. the price and the size of the forest, which are expected.
\begin{figure}[H]
\includegraphics[width=14cm]{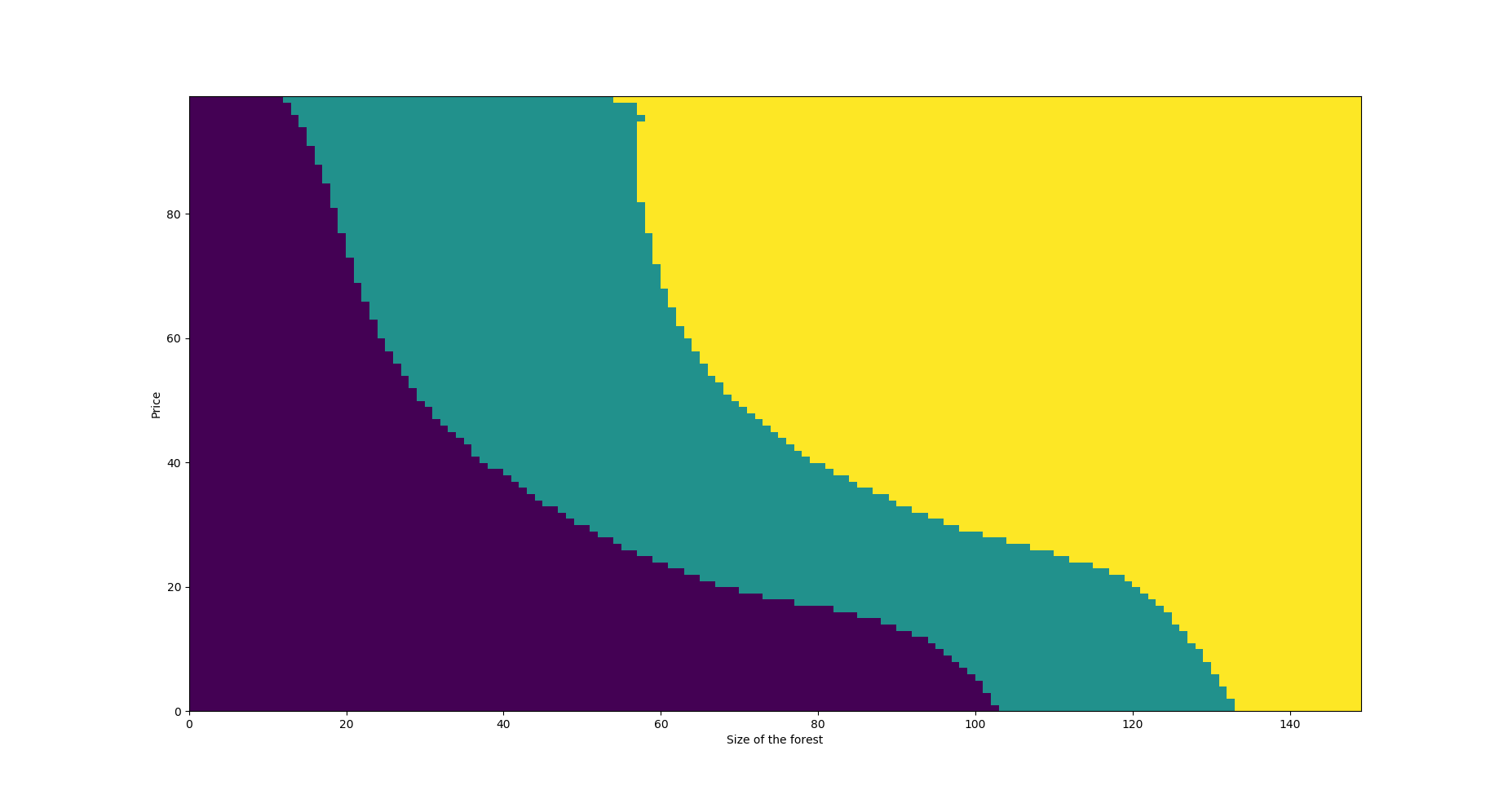}
\caption{The optimal strategy with respect to the price $P$ and the size of the forest $R$. The blue region corresponds to the plantation region, the yellow region corresponds to the harvesting region, the green region corresponds to the continuation region}\label{fig1}
\end{figure}
We note that the region to harvest is increasing with the price, and the region to renew is decreasing with the price. We never plant and harvest in the same time.

We now study the sensitivity w.r.t. the different parameters. For that we will change parameter by parameter.

\begin{figure}[H]
\includegraphics[width=14cm]{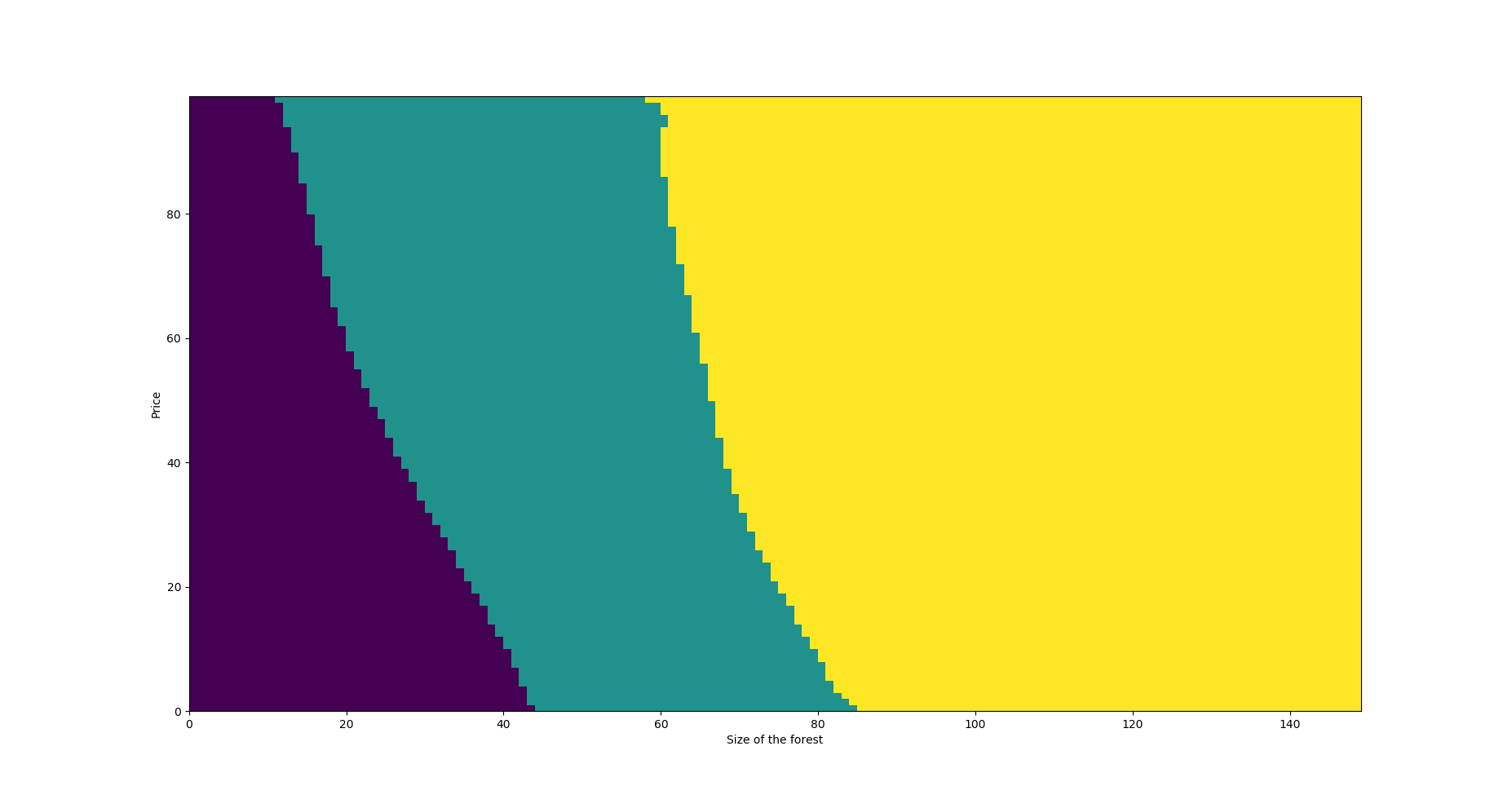}
\caption{In this figure the parameter $\lambda$ is now $0.9$}\label{fig2}
\end{figure}
If $\lambda$ is bigger in this case the region to harvest is more important and the region to renew is less important, since the growth is more important.

\begin{figure}[H]
\includegraphics[width=14cm]{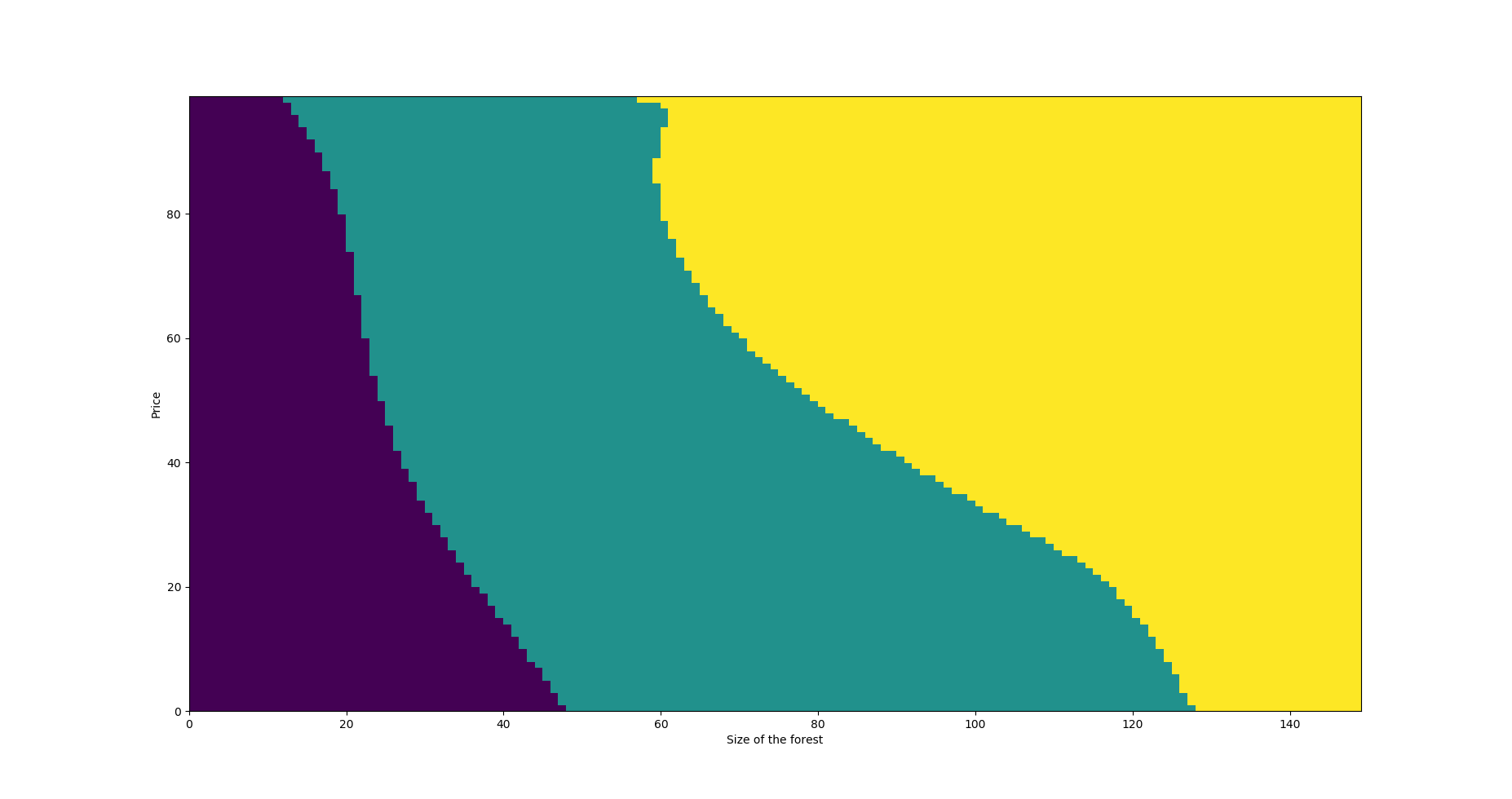}
\caption{In this figure the parameter $\eta$ is now $0.8$}\label{fig3}
\end{figure}
If $\eta$ is bigger in this case the region to renew is less important if the price is cheap, since the growth is slow and it is not interesting to renew except if the size of the forest is really small.

\begin{figure}[H]
\includegraphics[width=14cm]{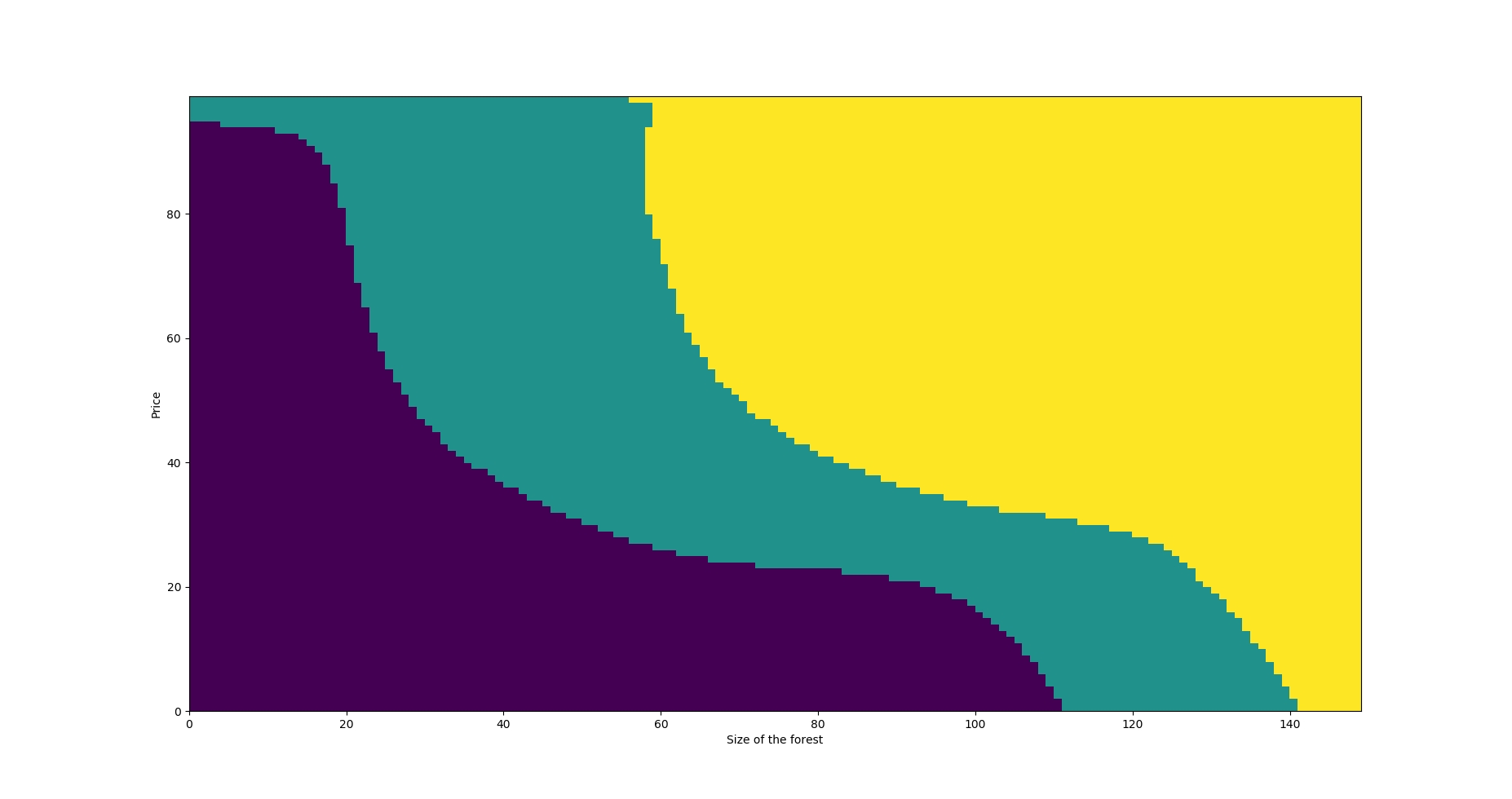}
\caption{In this figure the drift $\mu$ of the price is now $0.09$}\label{fig4}
\end{figure}
If the drift of the price is more important, the region to harvest is less important for a low price since the manager prefer to wait except if the size is too important because in this case the growth is negative, and the region to renew is more important because we know that the price will be better in the future.
\begin{figure}[H]
\includegraphics[width=14cm]{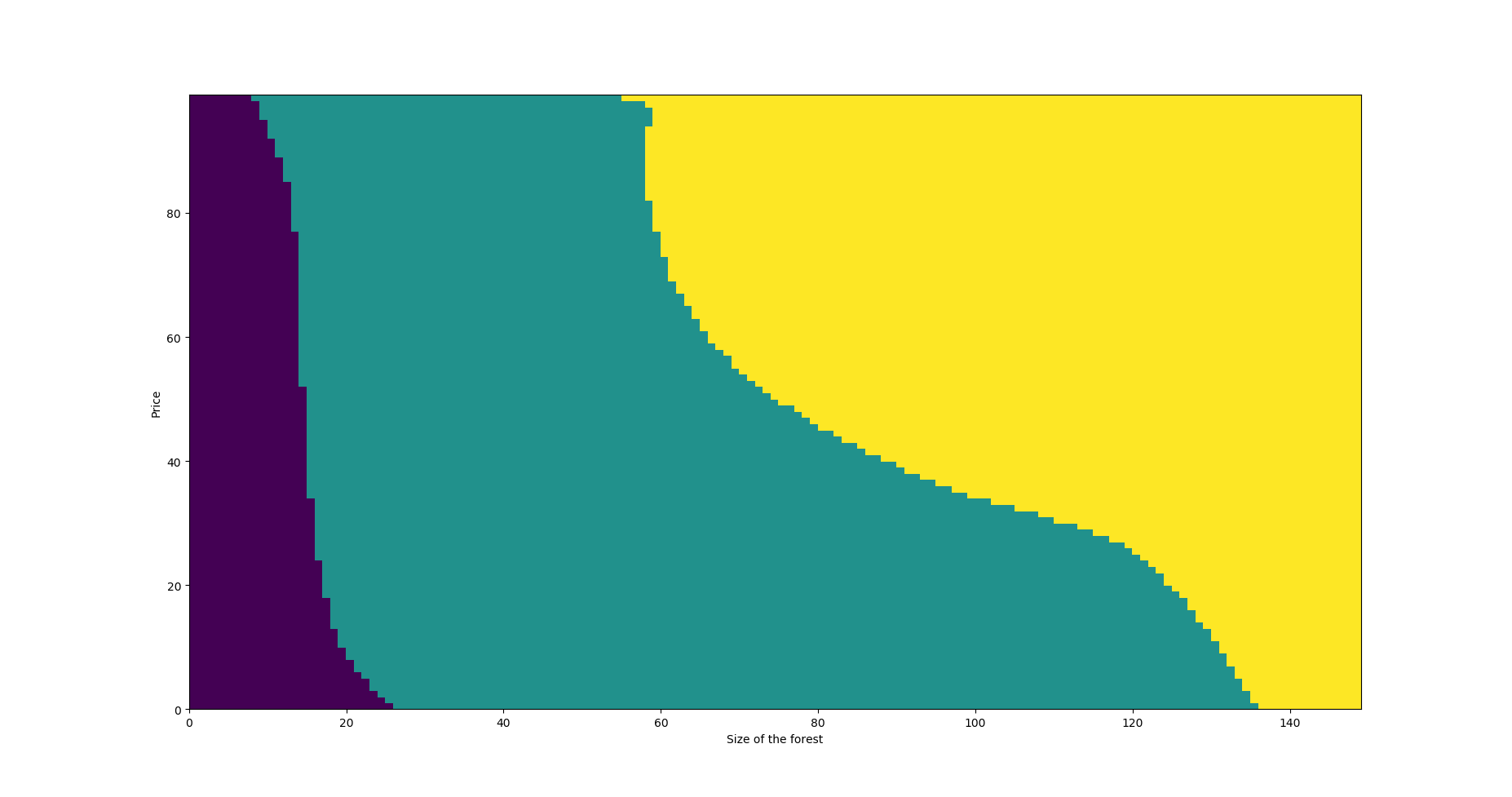}
\caption{In this figure the proportional costs $c_1$ and $c_3$ are now $0.15$}\label{fig5}
\end{figure}
If the costs are more expensive, the region to renew is less important because it es expensive to renew and harvest so we renew only if the size is really small.
\section{Proof of the main results}\label{sec5}

\subsection{Growth condition on $v$}\label{preuve1}
We provide in this subsection an upper-bound for the growth of the function $v$.

\vspace{2mm}

For any $(t,r)\in[0,T]\times\R_+$, we define the process  $\bar R^{t,r}$ by $\bar R^{t,r}_t=r$ and
 \beqs
 d \bar R^{t,r}_s & = & \eta \bar R^{t,r}_s(\lambda-\bar R^{t,r}_s)ds + \gamma \bar R^{t,r}_s  dB_s \;, \quad \forall \; s \in [t,T] \setminus \{ t_i ~:~ N(t) +1 \leq i \leq n\}\;,\\
 \bar R^{t,r}_{t_i} &=&\bar R^{t,r}_{t_i^-} + M \;, \quad \mbox{ for } ~N(t) +1 \leq i \leq n\;,
 \enqs
 where $M:= \max_{\xi\in[0,K]} g(\xi)+g_0$.
 We remark that the process $ \bar R^{t,r}$ can be written under the following form
 \beqs
 \bar R^{t,r}_s & = & r+ \int_t^s \eta \bar R^{t,r}_u(\lambda-\bar R^{t,r}_u)du + \int_t^s \gamma \bar R^{t,r}_u  dB_u+  \big(N(s)-N(t)\big)M \;,
 \enqs
 for $s\in[t,T]$. That corresponds to never harvest and renew always the maximum.

 We then have the following estimate on the process $\bar R ^{t,r}$.

\begin{Lemma}\label{lem-estim bar R}
For any $\ell\geq1$, there exists a constant $C_\ell$ such that
\beq\label{estim 2 bar R}
\E\Big[\sup_{s\in[t,T]}\big|\bar R^{t,r}_s\big|^\ell\Big] & \leq & C_\ell\big(1+|r|^{\ell}\big)\;,
\enq
for all $(t,r)\in[0,T]\times \R_+$.
\end{Lemma}
\ni \textbf{Proof.}
We first prove that for any $\ell\geq 1$, there exists a constant $C_\ell$ such that
\beq\label{estim 1 bar R}
\sup_{s\in[t,T]}\E\Big[\big|\bar R^{t,r}_s\big|^\ell\Big] & \leq & C_\ell\big(1+|r|^{\ell}\big)\;,
\enq
for all $(t,r)\in[0,T]\times \R_+$.
We argue by induction and we prove that for each $i=N(t),\ldots,n-1$ there exists a constant $C_{\ell,i}$ such that
\beq\label{estim 1bis bar R}
\E\Big[\big|\bar R^{t,r}_s\big|^\ell\Big] & \leq & C_{\ell,i}\big(1+|r|^{\ell}\big)\;,
\enq
for all $r\in \R_+$ and $s\in[t_i\vee t, (t_{i+1}\vee t)\wedge T)$.

\vspace{2mm}

\ni $\bullet$ For $i=N(t)$, using the closed formula of the logistic diffusion, we have
\beqs
\bar R^{t,r}_s & = & \frac{e^{(\eta\lambda-{\gamma^2\over 2}) (s-t) + \gamma (B_s-B_{t})}}{{ 1\over r} + \eta \int_{t}^s e^{(\eta\lambda-{\gamma^2\over 2}) (u-t)  + \gamma (B_u-B_{t})} du } \;,
\enqs
for all $s\in[ t, t_{N(t)+1}\wedge T)$.
Therefore we get
\beqs
\E\Big[\big|\bar R^{t,r}_s \big|^\ell\Big] & \leq & |r |^\ell\E\Big[\big|e^{(\eta\lambda-{\gamma^2\over 2}) (s-t) + \gamma (B_s-B_{t})}\big|^\ell\Big]\\
 & \leq &  |r|^\ell e^{(\ell|\eta\lambda-{\gamma^2\over 2}|+{|\ell\gamma|^2\over 2})(T-t)}
\enqs
for all $s\in[ t, t_{N(t)+1} \wedge T)$. Therefore \reff{estim 1bis bar R} holds true.

\vspace{2mm}

\ni $\bullet$ Suppose that the property holds for $i-1$. Still using the closed formula of the logistic diffusion, we have
\beqs
\bar R^{t,r}_s & = & \frac{e^{(\eta\lambda-{\gamma^2\over 2})  (s-t_i) + \gamma (B_s-B_{t_i})}}{{ 1\over \bar R^{t,r}_{t_i^-}+M} + \eta \int_{t_i}^s e^{(\eta\lambda-{\gamma^2\over 2}) (u-t_i)  + \gamma (B_u-B_{t_i})} du }\;,
\enqs
for all $s\in[t_i\vee t, (t_{i+1}\vee t)\wedge T)$.
Therefore we get
\beqs
\E\Big[\big|\bar R^{t,r}_s \big|^\ell\Big] & \leq & \E\Big[\big|(\bar R^{t,r}_{t_i^-}+M)e^{(\eta\lambda-{\gamma^2\over 2}) (s-t_i) + \gamma(B_s-B_{t_i})}\big|^\ell\Big]\\
 & \leq & \E\Big[\big|\bar R^{t,r}_{t_i^-}+M\big|^\ell\Big]e^{(\ell|\eta\lambda-{\gamma^2\over 2}|+{|\ell\gamma|^2\over 2})(T-t_{i})}\\
  & \leq & C'\Big( 1+ \E\Big[\big|\bar R^{t,r}_{t_i^-}\big|^\ell\Big]\Big)\;.
\enqs
Using the induction assumption and Fatou's Lemma, we get the result, and \reff{estim 1bis bar R} holds true for each $i=N(t),\ldots,n$. Taking $C_\ell=\max_{N(t)\leq i\leq n} C_{\ell,i}$, we get \reff{estim 1 bar R}.

\vspace{2mm}

We now prove \reff{estim 2 bar R}. Still using the closed formula of the logistic diffusion we have
\beqs
\big|\bar R^{t,r}_s \big|^\ell & \leq & \max_{N(t)\leq i\leq n} \big|(\bar R^{t,r}_{t_i^-}+M)\sup_{u\in[t_i\vee t, (t_{i+1}\vee t)\wedge T)}e^{(\eta\lambda-{\gamma^2\over 2})  (u-t_i) + \gamma (B_u-B_{t_i})}\big|^\ell\\
 & \leq & \sum_{i=N(t)}^{n}\big|(\bar R^{t,r}_{t_i^-}+M)\sup_{u\in[t_i\vee t, (t_{i+1}\vee t)\wedge T)}e^{(\eta\lambda-{\gamma^2\over 2}) (u-t_i) + \gamma (B_u-B_{t_i})}\big|^\ell\;,
\enqs
for all $s\in[t,T]$. Therefore, we get from the independence of $(B_u-B_{t_i})_{u\geq t_i}$ with $\Fc_{ t_i}$ and \reff{estim 1 bar R}
\beqs
\E\Big[\sup_{s\in[t,T]}\big|\bar R^{t,r}_s\big|^\ell\Big] & \leq & C \Big[\sum_{i=N(t)+1}^{n}\E\Big[\big|\bar R^{t,r}_{t_i^-}+M\big|^\ell\Big] + (1+|r|^\ell)\Big]
\\
 & \leq & C'_\ell ( 1+ |r|^\ell)\;,
\enqs
for some constant $C'_\ell$.
\ep

 \begin{Proposition}\label{estim R somme zeta}
 (i) For any $\ell\geq 1$, there exists a constant $C_\ell$ such that
 \beqs
 \E\Big[\sup_{s\in[t,T]}\big|R_s ^{t,r,\alpha}\big|^\ell\Big] & \leq & C_\ell \big(1+ |r|^\ell\big)
 \enqs
  for any strategy $\alpha \in \hat \Ac_{t,z,d}$.

\ni (ii) There exists a constant $C$ such that
\beqs
 \E\Big[ \Big(\sum_{k\geq 1} \zeta_k\mathds{1}_{\tau_k\leq T}\Big)^2\Big] & \leq & C \big(1+ |r|^4\big)
 \enqs
  for any strategy $\alpha \in \hat \Ac_{t,z,d}$.

 \end{Proposition}
\ni \textbf{Proof.} (i) Fix $\alpha=(t_i,\xi_i)_{N(t-\delta) +1 \leq i \leq n} \cup (\tau_k,\zeta_k)_{k \geq 1} \in \hat \Ac_{t,z,d}$. Using the definition of $\bar R ^{t,r}$  we have
\beqs
0~~\leq~~ R^{t,r,\alpha}_s & \leq & \bar R^{t,r}_s
\enqs
for all $s\in[t,T]$. Therefore we get from Lemma \ref{lem-estim bar R}
\beqs
 \E\Big[\sup_{s\in[t,T]}\big|R ^{t,r,\alpha}_s\big|^\ell\Big] & \leq & \E\Big[\sup_{s\in[t,T]} \big|\bar R^{t,r}_s\big|^\ell\Big] ~~  \leq ~~  C_\ell \big(1+ |r|^\ell\big)\;.
\enqs
(ii) We turn to the second estimate. From the dynamics \reff{eq croissance 1 bis} of $R^{t,r,\alpha}$, and since $R_T^{t,r,\alpha}\geq 0$ we have
\beqs
 \sum_{k\geq 1} \zeta_k\mathds{1}_{\tau_k\leq T} & \leq & r+ \int_t^T \eta R^{t,r,\alpha}_u(\lambda-R^{t,r,\alpha}_u)du + \int_t^T \gamma R^{t,r,\alpha}_u  dB_u + n M \;,
\enqs
where we recall that  $M=\max_{\xi\in[0,K]}g(\xi)+g_0$. Therefore, we get
\beqs
\E\Big[\Big(\sum_{k\geq 1} \zeta_k\mathds{1}_{\tau_k\leq T}\Big)^2\Big] & \leq & 4\Big(|r|^2+ \E\Big[\Big|\int_t^T \eta R^{t,r,\alpha}_u(\lambda-R^{t,r,\alpha}_u)du\Big|^2 + \Big|\int_t^T \gamma R^{t,r,\alpha}_u  dB_u\Big|^2 \Big]+ n^2 M^2\Big)\;.
\enqs
Therefore there exists a constant $C$ depending only on $T$, $\eta$, $\lambda$, $\gamma$, $M$ and $n$ such that
\beqs
\E\Big[\Big(\sum_{k\geq 1} \zeta_k\mathds{1}_{\tau_k\leq T}\Big)^2\Big] &  \leq & C\Big(|r|^2+ 1+ \E\Big[\sup_{s\in [t,T]}|R^{t,r,\alpha}_s|^4\Big]\Big) \;.
\enqs
Using estimate (i) we get the result.
\ep

\vspace{2mm}

We turn to the proof of the growth estimation for the value function $v$.

\vspace{2mm}

 \ni\textbf{Proof of Proposition \ref{Bound value function}.}  Fix $(t,z,d)\in \Dc$. From the definition of the function $L$ and the dynamics \reff{dyn X} and \reff{dyn P} of $X$ and $P$ we have
 \beqs
 \E\Big[ L\big(Z^{t,z,\alpha}_T\big)\Big] & \leq  &  \E\Big[X^{t,z,\alpha}_T\Big]+  \E\Big[\big|P^{t,p}_T\big|^2\Big]+  \E\Big[\big|R^{t,r,\alpha}_T\big|^2\Big]\\
  & \leq & x+ \E\Big[\sup_{s\in[t,T]} |P^{t,p}_s|^2 \Big] + \E\Big[\Big(\sum_{k\geq 1}\zeta_k\mathds{1}_{t\leq \tau_k\leq T}\Big)^2\Big] +  \E\Big[\big|R^{t,r,\alpha}_T\big|^2\Big]\\
   & & + e^{(2\mu+\sigma^2)(T-t)}|p|^2
\enqs
for any strategy $\alpha=(t_i,\xi_i)_{N(t-\delta) +1 \leq i \leq n} \cup (\tau_k,\zeta_k)_{k \geq 1} \in \hat \Ac_{t,z,d}$. From classical estimates there exists a constant $C$ such that
\beqs
\E\Big[\sup_{s\in[t,T]} |P^{t,p}_s|^2 \Big] & \leq & C\Big(1+ |p|^2\Big)
\enqs
for all $p\in\R_+^*$.  Using this estimate and Proposition \ref{estim R somme zeta}  we get
\beqs
v(t,z,d) & \leq & x+C\big(1+|r|^4+|p|^4+|q|^4\big)\;.
\enqs
Then by considering the strategy $\alpha^0=d\in \hat \Ac_{t,z,d}$ with no more intervention than $d$, we get
 \beqs
x 
& \leq & J(t,z,\alpha^0)~~\leq~~v(t,z,d)\;.
 \enqs
\ep

\subsection{Dynamic programming principle}\label{preuve2}
Before proving the dynamic programming principle, we need the
following results.
\begin{Lemma}\label{rem-ppte-model}{\rm
For any $(t,z,d)\in\Dc$ and any control $\alpha\in\hat \Ac_{t,z,d}$ we have the following properties.
\begin{enumerate}[(i)]
\item The pair $(Z^{t,z,\alpha}, d(.,\alpha))$ satisfies the following Markov property
\beqs
\E\big[ \phi(Z^{t,z,\alpha}_{\vartheta_2}) \big| \Fc_{\vartheta_1} \big] &=& \E\big[ \phi(Z^{t,z,\alpha}_{\vartheta_2}) \big| (Z^{t,z,\alpha}_{\vartheta_1},d(\vartheta_1,\alpha)) \big]
 \enqs
 for any bounded measurable function $\phi$, and any $\vartheta_1,\vartheta_2\in\Tc_{[t,T]}$ such that $\P\big(\vartheta_1\leq\vartheta_2\big)=1$.

 \item Causality of the control
 \beqs
 \alpha^\vartheta \in \hat \Ac_{\vartheta,Z_\vartheta^{t,z,d}, d(\vartheta,\alpha)} \quad & \text{ and }  & \quad d(\vartheta,\alpha) \in D_\vartheta \quad a.s.
 \enqs
for any $\vartheta\in\Tc_{[t,T]}$ where we set $\alpha^\vartheta =(t_{i},\xi_{i})_{  N(\vartheta-\delta)+1\leq i\leq n}\cup(\tau_{k},\zeta_{k})_{k\geq \kappa(\vartheta,\alpha)+1}$ and
 \beqs
\kappa(\vartheta,\alpha) & = & \#\big\{ k\geq 1 ~:~\tau_k< \vartheta\big\}\;.
 \enqs
 \item The state process $Z^{t,z,\alpha}$ satisfies the following flow property
\beqs
Z^{t,z,\alpha} & = & Z^{\vartheta, Z^{t,z,\alpha}_{\vartheta},\alpha^{\vartheta}} \quad \text{on } [\vartheta,T]\;,
\enqs
for any $\vartheta\in\Tc_{[t,T]}$.
 \end{enumerate}}
\end{Lemma}
\ni \textbf{Proof.} These properties are direct consequences of the dynamics of $Z^{t,z,\alpha}$. \ep

\vspace{2mm}
We turn to the proof of the dynamic programming principles (DP1) and (DP2).
Unfortunately, we have not enough information on the value function $v$ to directly prove these results. In particular, we do not know the measurability of $v$ and this prevents us from computing expectations involving $v$ as in (DP1) and (DP2). We therefore provide weaker dynamic programing principles involving the envelopes $v_*$ and $v^*$ as in \cite {BT11}. Since we get the continuity of $v$ at the end, these results implies (DP1) and (DP2).
\begin{Proposition}\label{PropDP1}
For any $(t,z,d)\in\Dc$ we have
\beqs
v(t,z,d)& \geq & \sup_{\alpha \in \hat \Ac_{t,z,d}} \sup_{\vartheta \in \Tc_{[t,T]}} \E \Big[   v_*(\vartheta,Z_\vartheta^{t,z,d},d(\vartheta,\alpha))\Big]\;.
\enqs
\end{Proposition}

\ni \textbf{Proof.}
Fix $(t,z,d) \in \Dc$,  $\alpha \in \hat \Ac_{t,z,d}$ and  $\vartheta \in \Tc_{[t,T]}$. By definition of the value function $v$, for any $\varepsilon >0$ and $\omega \in \Omega$, there exists $\alpha^{\varepsilon, \omega} \in \hat \Ac_{\vartheta(\omega),Z^{t,z,\alpha}_{\vartheta(\omega)}(\omega),d(\vartheta(\omega), \alpha)}$, which is an $\varepsilon$-optimal control at $(\vartheta,Z_\vartheta^{t,z,\alpha}, d(\vartheta, \alpha))(\omega)$, i.e.
\beqs
v\big(\vartheta(\omega),Z^{t,z,\alpha}_{\vartheta(\omega)}(\omega), d(\vartheta(\omega), \alpha(\omega))\big) - \varepsilon &\leq &
J(\vartheta(\omega),Z^{t,z,\alpha}_{\vartheta(\omega)}(\omega), \alpha^{\varepsilon,\omega})\;.
\enqs
By a measurable selection theorem (see e.g. Theorem 82 in the appendix of Chapter III in \cite{delmey75}) there exists $\bar \alpha_\varepsilon =(t_i,\bar \xi_i)_{N(\vartheta)+1\leq i\leq n}\cup(\bar \tau_k,\bar \zeta_k)_{k\geq 1}\in \hat\Ac_{\vartheta,Z^{t,z,\alpha}_{\vartheta},d(\vartheta,\alpha)}$ s.t. $\bar \alpha_{\varepsilon}(\omega)= \alpha_{\varepsilon, \omega}(\omega)$ a.s., and so
\beq \label{inegalite DP1 omega}
v\big(\vartheta,Z^{t,z,\alpha}_{\vartheta}, d(\vartheta, \alpha)\big) - \varepsilon &\leq & J(\vartheta,Z^{t,z,\alpha}_{\vartheta}, \bar \alpha_{\varepsilon})\;,\qquad \P-a.s.
\enq
We now define by concatenation the control strategy $\bar \alpha$ consisting of the impulse control components of $\alpha$ on $[t,\vartheta)$, and the impulse control components $\bar \alpha_\varepsilon$ on $[\vartheta,T]$. By construction of the control $\bar \alpha$
 we have  $\bar \alpha\in \hat\Ac_{t,z,d}$, $Z^{t,z,\bar \alpha} =   Z^{t,z,\alpha}$ on $[t,\vartheta)$,
  $d(\vartheta,\bar \alpha)= d(\vartheta,\alpha)$, and $\bar \alpha^\vartheta=\bar \alpha_\varepsilon$.
  From Markov property, flow property, and causality features of our model, given by Lemma \ref{rem-ppte-model},  the definition of the performance criterion and the law of iterated conditional expectations,  we get
\beqs
J(t,z,\bar \alpha)&=& \E \Big[   J(\vartheta,Z_\vartheta^{t,z,\alpha}, \bar \alpha_\varepsilon)\Big]\;.
\enqs
Together with \reff{inegalite DP1 omega}, this implies
\beqs
v(t,z,d) & \geq & J(t,z,\bar \alpha)\\
 &  \geq & \E \Big[  v_{*}(\vartheta,Z^{t,z,\alpha}_{\vartheta},d(\vartheta,\alpha))\Big] - \varepsilon \;.
\enqs
Since $\eps$, $\vartheta$ and $\alpha$ are arbitrarily chosen, we get the result.
\ep

\vspace{2mm}

We now prove (DP2), which is equivalent to the following proposition.
\begin{Proposition}\label{DP2}
For all  $(t,z,d) \in \Dc$, we have
\beqs
v(t,z,d)& \leq & \sup_{\alpha \in \hat \Ac_{t,z,d}} \inf_{\vartheta \in \Tc_{[t,T]}} \E \Big[   v^*\big(\vartheta,Z_\vartheta^{t,z,\alpha},d(\vartheta,\alpha)\big)\Big]\;.
\enqs
\end{Proposition}

\ni \textbf{Proof.}
Fix $(t,z,d) \in \Dc$, $\alpha \in \hat \Ac_{t,z,d}$ and $\vartheta \in \Tc_{[t,T]}$.  From the definitions of the performance criterion and the value functions, the law of iterated conditional expectations, Markov property, flow property, and causality features  of our model given by Lemma \ref{rem-ppte-model}, we get
\beqs
J(t,z,\alpha) & = &  \E\Big[ \E\Big[L \Big( Z^{\vartheta,Z^{t,z,\alpha}_\vartheta,\alpha^\vartheta}_T\Big) \Big| \Fc_\vartheta \Big]\Big]
~~= ~~  \E \Big[ J\big(\vartheta, Z_\vartheta^{t,z,\alpha},\alpha^\vartheta\big)\Big]\\
&\leq&  \E \Big[  v^*\big(\vartheta,Z_\vartheta^{t,z,\alpha},d(\vartheta,\alpha)\big)\Big]\;.
\enqs
Since $\vartheta$ and $\alpha$ are arbitrary, we obtain the required inequality.
\ep
\subsection{Viscosity properties}\label{preuve3}

We first need the following comparison result. We recall that $\Zc=\R\times\R_+\times\R_+^*\times\R_+^*$ and  $D_{t}$ is given by \reff{def D_t}.

\begin{Proposition}\label{Lemma-comp}
Fix $k\in\{0,\ldots,m-1\}$ (resp. $k\in\{m,\ldots,n-1\}$) and $g:\Zc\times D_{t_{k+1}}\rightarrow\R$ a continuous function. Let $\underline w:\Dc_k\rightarrow \R$ a viscosity  subsolution to \reff{EDP0}-\reff{EDP1} and
\beq\label{CondTerm1}
\underline w(t_{k+1}^-,z,d) & \geq & \max\big\{\Nc_2 g(z,d)~,~\bar \Nc_2 g(z,d)\big\} \;,\quad (z,d)\in\Zc\times D_{t_{k+1}}\\ \nonumber
\mbox{( resp. } \underline w(t_{k+1}^-,z,d) & \geq & \max\big\{\Nc_1 g(z,d)~,~\bar \Nc_1 g(z,d)\big\} \;,\quad (z,d)\in\Zc\times D_{t_{k+1}} \mbox{ ) } \;,
\enq
and $\bar w:\Dc_k\rightarrow \R$ a viscosity supersolution to \reff{EDP0}-\reff{EDP1}
\beq\label{CondTerm2}
\bar w(t_{k+1}^-,z,d) & \leq & \max\big\{\Nc_2 g(z,d)~,~\bar \Nc_2 g(z,d)\big\} \;,\quad (z,d)\in\Zc\times D_{t_{k+1}}\\ \nonumber
\mbox{( resp. } \bar w(t_{k+1}^-,z,d) & \leq & \max\big\{\Nc_1 g(z,d)~,~\bar \Nc_1 g(z,d)\big\} \;,\quad (z,d)\in\Zc\times D_{t_{k+1}} \mbox{ ) } \;.
\enq
Suppose there exists a constant $C>0$ such that
\beq\label{growth-cond1}
\underline w(t,z,d) & \leq & x+C\big(1+|r|^4+|p|^4+|q|^4+|d|^4\big)\\ \label{growth-cond2}
\bar w(t,z,d) & \geq & x
\;,
\enq
for all $(t,z,d)\in\Dc_k$ with $z=(x,r,p,q)$.
 Then $\underline w\leq \bar w$ on  $\Dc_k$.
In particular there exists at most a unique viscosity solution  $w$ to \reff{EDP0}-\reff{EDP1}-\reff{CondTerm1}-\reff{CondTerm2},  satisfying \reff{growth-cond1}-\reff{growth-cond2} and $w$ is continuous on $[t_k,t_{k+1})\times  \Zc$.
\end{Proposition}
The proof is postponed to the end of this section.
We are now able to state viscosity properties and uniqueness of $v$.

\paragraph{Viscosity property on $\Dc^1_k$.} Fix $k=0,\ldots,n-1$ and $(t,z,d)\in\Dc^1_k$ with $z=(x,r,p,q)$ and $r=0$.

1) We first prove the viscosity supersolution. Let $\varphi\in C^{1,2}(\Dc_k)$ such that
\beq\label{cond-fct-test r=0sursol}
(v_*-\varphi)(t,z,d) & = & \min_{\Dc_k} (v_*-\varphi) \;.
\enq
Consider a sequence $(s_\ell,z_\ell,d_\ell)_{\ell \in \N}$ of $\Dc_k$ such that
\beqs
\big(s_\ell,z_\ell,d_\ell,v(t_\ell,z_\ell,d_\ell)\big) & \xrightarrow[\ell\rightarrow+\infty]{} & \big(t,z,d,v_*(t,z,d)\big) \;.
\enqs
Applying Proposition \ref{PropDP1} with $\vartheta=s_\ell+h_\ell$ where $h_\ell \in(0,s_{\ell+1}-s_\ell)$. We have for $\ell$ large enough
\beqs
v(s_\ell,z_\ell,d_\ell) & \geq & \E\Big[v_*(s_\ell+h_\ell,Z^\ell_{s_\ell+h},d_\ell)\Big]\;,
\enqs
where $Z^\ell$ stands for $Z^{s_\ell,z_\ell,\alpha^0}$ with $\alpha^0$ the strategy with no more interventions than $d$. From \reff{cond-fct-test r=0sursol}, we get
\beqs
\chi_\ell +\varphi(s_\ell,z_\ell,d_\ell)& \geq & \E\Big[\varphi(s_\ell+h_\ell,Z^\ell_{s_\ell+h_\ell},d_\ell)\Big]\;,
\enqs
with $\chi_\ell := v(s_\ell,z_\ell,d_\ell)-v_*(t,z,d) - \varphi(s_\ell,z_\ell,d_\ell) +\varphi(t,z,d)\rightarrow0$ as $\ell\rightarrow \infty$.
Taking $h_\ell=\sqrt{|\chi_\ell|}$ and applying Ito's formula we get
\beqs
{1\over h_\ell} \E\Big[\int_{s_\ell}^{s_\ell+h_\ell}-\Lc\varphi(s,Z^\ell_{s},d_\ell)ds\Big] & \geq & -\sqrt{|\chi_\ell|}\;.
\enqs
Sending $\ell$ to $\infty$, we get the supersolution property from the mean value theorem.

\vspace{2mm}

2) We turn to the viscosity subsolution. Let $\varphi\in C^{1,2}(\Dc_k)$ such that
\beq\label{cond-fct-test r=0soussol}
(v^*-\varphi)(t,z,d) & = & \max_{\Dc_k} (v^*-\varphi) \;.
\enq
Consider a sequence $(s_\ell,z_\ell,d_\ell)_{\ell \in \N}$ of $\Dc_k$ such that
\beqs
\big(s_\ell,z_\ell,d_\ell,v(s_\ell,z_\ell,d_\ell)\big) & \xrightarrow[\ell\rightarrow+\infty]{} & \big(t,z,d,v^*(t,z,d)\big) \;.
\enqs
From Proposition \ref{DP2} we can find for each $\ell \in \N$ a control $\alpha^\ell=(t_i,\xi^\ell_i)_{N(t_\ell-\delta)+1\leq i\leq n}\cup (\tau_k,\zeta_k)_{k\geq1}\in\hat \Ac_{s_\ell,z_\ell,d_\ell}$ such that
\beqs
v(s_\ell,z_\ell,d_\ell) & \leq & \E\Big[v^*(s_\ell+h_\ell,Z^\ell_{s_\ell+h_\ell},d)\Big]+{1\over \ell} \;,
\enqs
where $Z^\ell$ stands for $Z^{s_\ell,z_\ell,\alpha^\ell}$ and $h_\ell\in (0,s_{\ell+1}-t_\ell)$ is a constant that will be chosen later.

We first notice that
\beq\label{conv zero r alpha}
\sup_{s\in[s_\ell,s_\ell+h_\ell]}|R^\ell_s| & \xrightarrow[\ell\rightarrow \infty]{\P-a.s.} & 0 \;.
\enq
Indeed, we have
\beq\label{majRell}
0~~ \leq ~~ R^\ell_s &  \leq &  \bar R^\ell_s\;,\quad s\geq s_\ell
\enq
 where $\bar R^\ell$ is given by
\beqs
\bar R^\ell_s & = & r_\ell+\int_{s_\ell}^s \eta\bar R^\ell_u(\lambda-\bar R^\ell_u)du+\int_{s_\ell}^s\bar R^\ell_udB_u\;,\quad s\geq s_\ell \;.
\enqs
Since $r_\ell\xrightarrow[\ell\rightarrow \infty]{} r$ (and $r=0$), we have $\sup_{s\in[s_\ell,s_\ell+h_\ell]}|\bar R^\ell_s| \xrightarrow[\ell\rightarrow \infty]{}0$ as $\ell\rightarrow \infty$ and we get \reff{conv zero r alpha}. In particular, we deduce that up to a subsequence
\beq\label{conv somme zeta}
\sum_{k\geq1} \zeta_k^\ell\mathds{1}_{\tau_k^\ell\leq s_\ell+h_\ell} & \xrightarrow[\ell\rightarrow+\infty]{\P-a.s.} & 0\;.
\enq
Indeed, we have from \reff{eq croissance 1 bis} and \reff{majRell}
\beq
\sum_{k\geq1} \zeta_k^\ell\mathds{1}_{\tau_k^\ell\leq s_\ell+h_\ell} & \leq & r_\ell+\int_{s_\ell}^{s_\ell+h_\ell}\eta \lambda R^\ell_udu+\int_{s_\ell}^{s_\ell+h_\ell}\eta R^\ell_udB_u\nonumber\\\label{majsomme coupe}
 & \leq & r_\ell+h_\ell\eta\lambda\sup_{s\in[s_\ell,s_\ell+h_\ell]}|\bar R^\ell_s| +\big|\int_{s_\ell}^{s_\ell+h_\ell}\eta R^\ell_udB_u\big|\;.\nonumber
\enq
From BDG inequality and \reff{majRell}, we get from \reff{conv zero r alpha}
\beqs
\E\Big[\big|\int_{s_\ell}^{s_\ell+h_\ell}\eta R^\ell_udB_u\big|\big] & \xrightarrow[\ell\rightarrow+\infty]{} & 0\;,
\enqs
and hence, up to a subsequence $\big|\int_{s_\ell}^{s_\ell+h_\ell}\eta R^\ell_udB_u\big|\rightarrow0$ as $\ell\rightarrow+\infty$.
From this convergence \reff{conv zero r alpha} and  \reff{majsomme coupe}, we get \reff{conv somme zeta}.

 We then define the process $\tilde X ^\ell$ by
\beqs
\tilde X ^\ell _s & = & x_\ell+\sum_{k\geq1} P_{\tau_k^\ell}\zeta_k^\ell\mathds{1}_{\tau_k^\ell\leq s}
\enqs
and observe that   from \reff{conv somme zeta}
\beq\label{Xtilde domine X}
\tilde X ^\ell_{s_\ell+h_\ell} & \xrightarrow[\ell\rightarrow+\infty]{\P-a.s.} & x\;,\\
\tilde X ^\ell & \geq & X ^\ell\;.\nonumber
\enq
Since $v$ is nondecreasing in the $x$ component, it is the same for $v^*$. We get
\beqs
v(s_\ell,z_\ell,d_\ell) & \leq & \E\Big[v^*(s_\ell+h_\ell,\tilde Z^\ell_{s_\ell+h_\ell},d)\Big]+{1\over \ell}
\enqs
where $\tilde Z^\ell=(\tilde X^\ell,R^\ell,P^\ell,Q^\ell)$. We then get
from \reff{cond-fct-test r=0soussol}
\beqs
\chi_\ell +\varphi(s_\ell,z_\ell,d_\ell)& \leq & \E\Big[\varphi(s_\ell+h,\tilde Z^\ell_{s_l+h},d_\ell)\Big] +{1\over \ell}\;,
\enqs
where
$\chi_\ell := v(s_\ell,z_\ell,d_\ell)-v^*(t,z,d) - \varphi(s_\ell,z_\ell,d_\ell) + \varphi(t,z,d)\rightarrow0$ as $\ell\rightarrow+\infty$.
Applying Ito's formula and taking $h_\ell=\sqrt{|\chi_\ell|}$ we get by sending $\ell$ to $\infty$ as previously
\beqs
-\Lc\varphi(t,z,d) & \leq & 0\;.
\enqs

\paragraph{Viscosity property on $\Dc^2_k$.} Fix $k=0,\ldots,n-1$ and $(t,z,d)\in\Dc^2_k$. Then $v(.,d)$ is the value function associated to an optimal impulse control problem with nonlocal operator $\Hc$. Using the same arguments as in the proof of Theorem 5.1 in \cite{MLP07}, we obtain that $v$ is a  viscosity solution to \reff{EDP1} on $\Dc_k^2$.

\paragraph{Viscosity property and continuity on $\{t_k\}\times \Zc\times D_{t_k}$.}

 We prove it by a backward induction on $k=0,\ldots,n$.

\ni $\bullet$ Suppose that $k=n$ $i.e.$ $t_k=T$.

1) We first prove the subsolution property. Fix some $z=(x,r,p,q)\in \Zc$ and $d=(t_i,e_i)_{n-m+1\leq i \leq n}\in D_{t_{n}}$ and consider a sequence $(s_\ell,z_\ell,d_\ell)_{\ell \in \N}$ with $z_\ell=(x_\ell,r_\ell,p_\ell,q_\ell)$ and $d_\ell=(t_i,e_i^\ell)_{n-m+1\leq i \leq n}$ such that
\beqs
(s_\ell,z_\ell,d_\ell,v(s_\ell,z_\ell,d_\ell)) & \xrightarrow[\ell\rightarrow+\infty]{} & (T^-,z,d,v_*(T^-,z,d)) \;.
\enqs

By considering a strategy $\alpha^\ell\in \hat \Ac_{s_\ell,z_\ell,d_\ell}$ with a single renewing order $(T,e)$ with $e \leq K$
and the stopping time $\vartheta=T$, we get from the definition of $v$
\beqs
v(s_\ell,z_\ell,d_\ell) & \geq &
\E\Big[L\Big(\Gamma^{rn}_1 \big(\Gamma^{rn}_2(Z^{s_\ell,z_\ell,\alpha^\ell}_{T^-},e^\ell_{n-m+1}),e \big) \Big)\Big] \;.
\enqs
From the continuity of the functions $L$, $\Gamma^{rn}_1$ and $\Gamma^{rn}_2$, we get
\beqs
L\Big(\Gamma^{rn}_1 \big(\Gamma^{rn}_2(Z^{s_\ell,z_\ell,\alpha^\ell}_{T^-},e^\ell_{n-m+1}),e \big)  & \xrightarrow[\ell\rightarrow+\infty]{\P-a.s.} & L\Big(\Gamma^{rn}_1 \big(\Gamma^{rn}_2(z,e_{n-m+1}),e \big) \;.
\enqs
From Fatou's Lemma and since $e\leq K$ is arbitrarily chosen, we get by sending $\ell$ to $\infty$
\beq\label{estim v_* 1}
v_*(T^-,z,d) & \geq & \Nc_1 L(z,d)\;.
\enq
Fix now $a\in[0, r]$ and denote $a_\ell= \min\{a,r_\ell\}$. By considering a strategy $\alpha^\ell$ with an immediate harvesting order $(s_\ell,r_\ell)$ and a single renewing order $(T,e)$ and $\vartheta=T$, we get from the definition of $v$
 \beqs
v(s_\ell,z_\ell,d_\ell) & \geq &
\E\Big[L\Big(\Gamma^{rn}_1 \big(\Gamma^{rn}_2 \big(Z^{s_\ell,\Gamma^c(z_\ell,r_\ell),\alpha^\ell}_{T^-},e_{n-m+1} \big),e \big) \Big)\Big] \;.
\enqs
From the continuity of the functions $L$, $\Gamma^c$, $\Gamma^{rn}_1$ and $\Gamma^{rn}_2$, we get
\beqs
L\Big(\Gamma^{rn}_1 \big(\Gamma^{rn}_2 \big(Z^{s_\ell,\Gamma^c(z_\ell,r_\ell),\alpha^\ell}_{T^-},e_{n-m+1} \big),e \big) \Big)  & \xrightarrow[\ell\rightarrow+\infty]{\P-a.s.} & L\Big(\Gamma^{rn}_1 \big(\Gamma^{rn}_2 \big(\Gamma^c(z,r),e_{n-m+1} \big),e \big) \Big)\;.
\enqs
 From Fatou's Lemma and since $e\leq K$ and $a\in[0,r]$ are arbitrarily chosen, we get by sending $\ell$ to $\infty$
\beq\label{estim v_*2}
v_*(T^-,z,d) & \geq & \bar \Nc_1 L(z,d)\;.
\enq
From \reff{estim v_* 1} and \reff{estim v_*2}, we get the subsolution property at $(T^-,z,d)$.

\vspace{2mm}

2) We turn to the supersolution property.
We argue by contradiction and suppose that there exist $z=(x,r,p,q)\in \Zc$ and $d\in D_{t_{n}}$  such that
\beqs
v^*(T^-,z,d) & \geq & \max\Big\{\Nc_1L(z,d)~,~\bar \Nc_1L(z,d)\Big\} + 2\eps\;,
\enqs
with $\eps>0$.
We fix a sequence $(s_\ell,z_\ell,d_\ell)_{\ell \in \N}$ in $\Dc$ such that
\beq\label{approxlimsup}
(s_\ell,z_\ell,d_\ell,v(s_\ell,z_\ell,d_\ell)) & \xrightarrow[\ell\rightarrow+\infty]{} & (T^-,z,d,v^*(T^-,z,d))\;.
\enq
We then can find  $s>0$ and a sequence of smooth functions $(\varphi^h)_{h\geq 1}$ on  $[T-s,T]\times \Zc \times D_{t_{n}}$ such that $\varphi^h\downarrow v^*$ on $[T-s,T)\times \Zc \times D_{t_{n}}$, $\varphi^h\downarrow v^*(.^-,.,.)$ on $\{T\}\times \Zc \times D_{t_{n}}$ as $h\uparrow+\infty$ and
  \beq\label{cond1 varphi}
\varphi^h(t',z',d') & \geq & \max\Big\{\Nc_1L(z',d')~,~\bar \Nc_1L(z',d')\Big\} + \eps\;,
\enq
on some neighborhood $\Bc^h$ of $(T,z,d)$ in $[t_{n},T]\times \Zc \times D_{t_{n}}$. Up to a subsequence, we can assume that $\Bc^h_\ell := [t_\ell,T]\times B((z_\ell,d_\ell), \delta_\ell^h)\subset \Bc^h$ for $\delta_\ell^h$ sufficiently small.
Since $v^*$ is locally bounded, there is some $\iota > 0$ such that $|v^*|\leq\iota$ on $\Bc^h$. We therefore get $\varphi^h\geq -\iota$ on $\Bc^h$. We then define the function $\varphi^h_\ell$ by
\beqs
\varphi^h_\ell(t',z',d') & = & \varphi^h(t',z',d') + 3 \iota {|(z',d')-(z_\ell,d_\ell)|^2\over |\delta^h_\ell|^2 }+\sqrt{T-t'} \;,
\enqs
and we observe that
\beq\label{cond2 varphi}
(v^*-\varphi_\ell^h) & \leq & -\iota~~<~~0 \quad \mbox{ on } [t_\ell,T]\times \partial B((z_\ell,d_\ell), \delta_\ell^h)\;.
\enq
Since ${\partial\sqrt{T-t}\over \partial t}\rightarrow -\infty$ as $t\rightarrow T^-$,  we can choose $h$ large enough such that
\beq\label{cond3 varphi}
-\Lc \varphi _\ell ^h & \geq & 0 ~\mbox{ on }~ \Bc_\ell^h\;.
\enq
From the definition of $v$ we can find  $\alpha^\ell =(t_i,\xi_i^\ell)_{N(t_\ell-\delta)+1\leq i\leq n} \cup (\tau_k^\ell,\zeta_k^\ell)_{k\geq 1}\in \hat \Ac_{s_\ell,z_\ell,d_\ell}$ such that
\beq\label{Dynpro1/n}
v(t_\ell,z_\ell,d_\ell) & \leq & \E\Big[ L\big(Z^{\ell}_{T} \big)\Big]+{1\over \ell}\;,
\enq
where $Z^\ell$ stands for $Z^{s_\ell,z_\ell,\alpha^\ell}$. Denote by $\theta_\ell^h  = \inf\{s\geq s_\ell~:~(s,Z^\ell,d_\ell)\notin \Bc_\ell^h\}\wedge \tau^\ell_1$.
From Ito's formula, \reff{cond1 varphi}, \reff{cond2 varphi} and \reff{cond3 varphi} we have
\beqs
\varphi_\ell^h(s_\ell,z_\ell,d_\ell) & \geq &
 \E\Big[\Big(v\big(T,\Gamma^{rn}(\Gamma^c(Z^\ell_{T^-}, \zeta_1^\ell),\xi_{n-m}^\ell),d_\ell\cup(t_{n-m},\xi^\ell_{n-m}) \big)\mathds{1}_{\tau_1^\ell = T}  \\
 & & \qquad+v^*\big(\theta^\ell_n,\Gamma^c(Z^\ell_{{\theta_\ell^h}^-}, \zeta_1^\ell),d_\ell) \big)\mathds{1}_{\tau_1^\ell < T}\Big)\mathds{1}_{\tau^\ell_1\leq \theta_\ell^h}\Big] \\
  & & +  \E\Big[\Big(v\big(T,\Gamma^{rn}(Z^\ell_{T^-}, \xi^\ell_{n-m}),d_\ell\cup(t_{n-m},\xi^\ell_{n-m}) \big)\mathds{1}_{ \theta_\ell^h=T}  \\
  & & \qquad+v^*\big(\theta_\ell^h,Z^\ell_{{\theta_\ell^h}^-},d_\ell \big)\mathds{1}_{ \theta_\ell^h<t_k}\Big)\mathds{1}_{\tau^1_\ell> \theta_\ell^h}\Big]
 +\eps\wedge\iota\;.
\enqs
From \reff{Dynpro1/n} and the Markov  property given by Lemma \ref{rem-ppte-model} (i), we get by taking the conditional expectation given $\Fc_{\theta_\ell^h}$,
\beqs
 v(t_\ell,z_\ell,d_\ell) & \leq &
 \E\Big[\Big(v\big(T,\Gamma^{rn}(\Gamma^c(Z^\ell_{T^-}, \zeta_1^\ell),\xi_k^\ell),d_\ell\cup(t_{n-m},\xi^\ell_{n-m}) \big)\mathds{1}_{\tau_1^\ell = T}\\
  & & \qquad +v^*\big(\theta_\ell^h,\Gamma^c(Z^\ell_{{\theta_\ell^h}^-}, \zeta_1^\ell),d_\ell) \big)\mathds{1}_{\tau_1^\ell < T}\Big)\mathds{1}_{\tau^\ell_1\leq \theta_\ell^h}\Big]  \\
& & +  \E\Big[\Big(v\big(T,\Gamma^{rn}(Z^\ell_{T^-}, \xi^\ell_{n-m}),d_\ell\cup(t_{n-m},\xi^\ell_{n-m}) \big)\mathds{1}_{ \theta_\ell^h=T}\\
 & & \qquad +v^*\big(\theta_\ell^h,Z^\ell_{{\theta_\ell^h}^-},d_\ell \big)\mathds{1}_{ \theta_\ell^h<T}\Big)\mathds{1}_{\tau_1^\ell> \theta_\ell^h}\Big]
 +{1\over \ell}  \;.
\enqs
We therefore get
\beqs
\varphi^h(s_\ell,z_\ell,d_\ell) +\sqrt{T-s_\ell}~~=~~\varphi_\ell^h(s_\ell,z_\ell,d_\ell) & \geq & v(s_\ell,z_\ell,d_\ell)+\eps\wedge\iota-{1\over \ell} \;.
\enqs
Sending $\ell$ and $h$ to $+\infty$ we get a contradiction with \reff{approxlimsup}.

 \vspace{2mm}

\ni $\bullet$ Suppose that the property holds true for $k+1$.   From Proposition \ref{Lemma-comp}, the function $v$ is continuous on $D_{t_{k+1}}$. Therefore, we get from Propositions \ref{PropDP1} and \ref{DP2}
\beqs
v(t,z,d) & = & \sup_{\alpha\in\hat \Ac_{t,z,d}}\E\Big[v\big(t_{k+1},Z^{t,z,\alpha}_{t_{k+1}},d(t_{k+1},\alpha)\big)\Big]
\enqs
for all $(t,z,d)\in\Dc_k$.

We can then apply the same arguments as for $k=n$ and we get the viscosity property at $(t_{k+1}^-,z,d)$ for all $(z,d)\in \Zc\times D_{t_{k+1}}$.

\paragraph{Proof of Proposition \ref{Lemma-comp}.}
We fix the functions $\underline w$ and $\bar w$ as in the statement of Proposition \ref{Lemma-comp}. We then introduce  as classically done a perturbation of $\bar w$ to make it a strict supersolution.

\begin{Lemma}
Consider the function $\psi$ defined by
\beqs
\psi(t,z,d) & = & x+pr+\tilde C_1e^{-\tilde C_2 t}\big( 1+|r|^4+|p|^4+|q|^4+|d|^4\big) \;,
\enqs
where $\tilde C_1$ and $\tilde C_2$ are two positive constants and define for $m\geq 1$ the function $\bar w_m$ on $\Dc_k$ by
\beqs
\bar w_m & = & \bar w +{1\over m}\psi\;.
\enqs
Then there exist $\tilde C_1$ and $\tilde C_2$ (large enough) such that the following properties hold.
  \begin{itemize}
\item The function $\bar w_m$ is a strict viscosity supersolution to \reff{EDP0}-\reff{EDP1} on $[t_k,t_{k+1})\times\Kc$ for any compact subset $\Kc$ of $\Zc\times D_{t_k}$ and any $m\geq 1$ : there exists a constant $\delta>0$ (depending on $\Kc$ and $m$) such that
\beqs
&&- \Lc \varphi(t,z,d) \geq  \delta\\
&&\mbox{(resp. }\min\Big\{ - \Lc \varphi(t,z,d) ~,~
\bar w_m(t,z,d) - \Hc \bar w_m(t,z,d)
\Big\} \geq  \delta\mbox{)}
\enqs
for any $(t,z,d)\in\Dc_k^1$ (resp. $(t,z,d)\in\Dc_k^2$) and $\varphi\in C^{1,2}(\Dc_k)$ such that $(z,d)\in \Kc$ and
\beqs
(\bar w_m-\varphi)(t,z,d) & = & \min_{\Dc_k} (\bar w_m-\varphi)\;.
\enqs
\item We have
\beq\label{lim-sursol-per}
\lim_{|(z,d)|\rightarrow+\infty}(\underline w-\bar w_m)(t,z,d) & = & -\infty\;.
\enq
\end{itemize}
\end{Lemma}
\ni \textbf{Proof.} A straightforward computation shows that
\beqs
\psi-\Hc\psi & \geq & c_2>0\;,
\enqs
on $\Dc_{k}$.
Since $\bar w$ is a viscosity supersolution to \reff{EDP1}, we get
\beq\label{estim obstacle wm}
\bar w_m-\Hc \bar w_m& \geq & {c_2\over m}~~=:~~\delta_0~~>~~0\;,
\enq
on $\Dc^2_k$.
Then, from the definition of  the operator $\Lc$ we get for $\tilde C _2$ large enough
\beqs
-\Lc\psi & > & 0\quad \mbox{ on } \Dc_{t_k}\;.
\enqs
In particular, since $-\Lc\psi$ is continuous, we get
\beq\label{estim -L psi}
\inf_{[t_k,t_{k+1})\times\Kc} -{1\over m}\Lc\psi~~=:~~Ê\delta_1 & > & 0
\enq
for any compact subset $\Kc$ of $\Zc\times D_{t_k}$.
By writing the viscosity supersolution property of $\bar w$, we deduce from \reff{estim obstacle wm} and \reff{estim -L psi} the desired strict viscosity supersolution property for $w_m$.

\ni Finally, from growth conditions \reff{growth-cond1} and \reff{growth-cond2}, we get \reff{lim-sursol-per} for $\tilde C_1$ large enough.
\ep

\vspace{3mm}

\ni To prove the comparison result, it suffices to prove that
\beqs
\sup_{\Dc_k} \;(\underline w - \bar w_m ) &  \leq & 0\;,
\enqs
for all $m\geq1$. We argue by contradiction and suppose that there exists $m\geq1$ such that
\beqs
\bar \Delta~~:=~~ \sup_{\Dc_k}\; (\underline w - \bar w_m ) &  > & 0\;.
\enqs
Since $\bar w_m-\underline w$ is u.s.c. on $\Dc_k$ and $\bar w_m-\underline w(t_{k+1}^-,.)\leq 0$, we get from  \reff{lim-sursol-per} the existence of an open subset $\Oc$ of $\Zc\times D_{t_k}$ and $(t_0,z_0,d_0)\in[t_k,t_{k+1})\times \Oc$ such that $\bar \Oc$ is compact and
\beqs
(\underline w-\bar w_m)(t_0,z_0,d_0) & = & \bar \Delta\;.
\enqs
We then consider the functions $\Phi_i$ and $\Theta_i$ defined on $[t_k,t_{k+1})\times \bar \Oc$ by
\beqs
\Phi_i(t,t',z,z',d,d') & = & \underline w(t,z,d)-\bar w_m(t',z',d')-\Theta_i(t,t',z,z',d,d')  \\
\Theta_i(t,t',z,z',d,d') & = & |t-t_0|^2+ |z-z_0|^{4}+|d-d_0|^{2}+    {i\over2}\big(|z-z'|^2+  |d-d'|^2\big)
\enqs
for all $(t,z,d),(t',z',d')\in D_k$ and $i\geq 1$. From the growth properties of $\underline w$ and $\bar w_m$, there exists $(\hat t_i,\hat t'_i,\hat z_i,\hat z'_i,\hat d_i,\hat d'_i)\in ([t_k,t_{k+1})\times \bar\Oc)^2$ such that
\beqs
\bar \Delta_i & := & \sup_{[t_k,t_{k+1})\times \bar\Oc}\Phi_i~~=~~\Phi_i(\hat t_i,\hat t'_i,\hat z_i,\hat z'_i,\hat d_i,\hat d'_i)\;.
\enqs
By classical arguments we get, up to a subsequence, the following convergences
\beq\nonumber
\big(\hat t_i,\hat t'_i,\hat z_i,\hat z'_i,\hat d_i,\hat d'_i,\big) & \xrightarrow[i\rightarrow+\infty]{} & \big(t_0,t_0,z_0,z_0,d_0,d_0, \big)\;,\\\nonumber
\Phi_i(\hat t_i,\hat t'_i,\hat z_i,\hat z'_i,\hat d_i,\hat d'_i) & \xrightarrow[i\rightarrow+\infty]{} &  (\underline w-\bar w_m)(t_0,z_0,d_0)\;,\\  \label{conv Theta i 0}
\Theta_i(\hat t_i,\hat t'_i,\hat z_i,\hat z'_i,\hat d_i,\hat d'_i)  & \xrightarrow[i\rightarrow+\infty]{} & 0\;.
\enq
In particular, we have $\max\{\hat t_i,\hat t_i'\}<T$ for $i$ large enough.
We then apply Ishii's Lemma (see Theorem 8.3 in \cite{CIL}) to $(\hat t_i,\hat t'_i,\hat z_i,\hat z'_i,\hat d_i,\hat d'_i)$ which realizes the maximum of  $\Phi_i$ and we get for any $\eps_i>0$, the existence of $(e_i,f_i,M_i)\in \bar J^{2,+}\underline w(\hat t_i,\hat z _i)$ and  $(e'_i,f'_i,M'_i)\in \bar J^{2,-}\bar w_m(\hat t'_i,\hat z' _i)$ such that
\beq\label{ishii-cond-1}
e_i ~~=~~{\partial \Theta_i\over \partial t}(\hat t_i,\hat t'_i,\hat z_i,\hat z'_i,\hat d_i,\hat d'_i) & & f_i~~=~~{\partial \Theta_i\over \partial z}(\hat t_i,\hat t'_i,\hat z_i,\hat z'_i,\hat d_i,\hat d'_i)\\ \label{ishii-cond-2}
e_i' ~~=~~{\partial \Theta_i\over \partial t'}(\hat t_i,\hat t'_i,\hat z_i,\hat z'_i,\hat d_i,\hat d'_i) & & f_i'~~=~~{\partial \Theta_i\over \partial z'}(\hat t_i,\hat t'_i,\hat z_i,\hat z'_i,\hat d_i,\hat d'_i)
\enq
and
\begin{equation}\label{ishii-cond-3}
\left(\begin{array}{cc}
M & 0 \\
0 & -M'
\end{array}\right)~~\leq~~Ê{\partial^2 \Theta_i\over \partial (z, z')^2}(\hat t_i,\hat t'_i,\hat z_i,\hat z'_i,\hat d_i,\hat d'_i)+{1\over i}\Big({\partial^2 \Theta_i\over \partial (z, z')^2}(\hat t_i,\hat t'_i,\hat z_i,\hat z'_i,\hat d_i,\hat d'_i)\Big)^2\;,
\end{equation}
for all $i\geq 1$. We then distinguish two cases.

\vspace{2mm}

\ni$\bullet$ Case 1: there exists a subsequence of $(\hat t_i,\hat t'_i,\hat z_i,\hat z'_i,\hat d_i,\hat d'_i)_{i \in \N}$ still denoted $(\hat t_i,\hat t'_i,\hat z_i,\hat z'_i,\hat d_i,\hat d'_i)_{i \in \N}$ such that
\beqs
(\hat t_i,\hat z_i,\hat d_i) & \in & \Dc_k^2\quad \mbox{ for all }~i\geq 1\;.
\enqs
From the viscosity subsolution property of $\underline w$ and the strict viscosity supersolution property of $\bar w_m$ we have
\beq\label{soussol}
\min\Big\{-\Lc[\hat z_i,\hat d_i,e_i,f_i,M_i]~;~(\underline w-\Hc \underline w)(\hat t_i,\hat z_i,\hat d_i)\Big\} & \leq & 0\\\label{surssol}
\min\Big\{-\Lc[\hat z'_i,\hat d'_i,e'_i,f_i',M_i']~;~(\bar w_m-\Hc \bar w_m)(\hat t_i,\hat z_i,\hat d_i)\Big\} & \geq & \delta\over m
\enq
where
\beqs
\Lc [z,d,e,f,M] & = & e + \mu p f_3+\rho qf_4+\eta r(\lambda-r)f_2 \\
 &  & +{1\over 2} \Big(\sigma^2p^2M_{3,3}+\varsigma^2q^2M_{4,4}+2\sigma\varsigma pq M_{3,4}+\gamma^2r^2 M_{2,2}\Big)
\enqs
for any $z\in \Zc$, $d\in D_{t_k}$, $e\in\R$, $f\in \R^4$  and any symmetric matrix $M\in \R^{4\times 4}$\;.
We then distinguish the following two possibilities in \reff{soussol}.

\vspace{2mm}

\ni \textbf{1.} Up to a subsequence we have
\beqs
w(\hat t_i ,\hat z _i,\hat d_i)-\Hc w(\hat t_i,\hat z_i,\hat d_i) \leq 0 ~~\mbox{  for all }~i\geq 1.
\enqs
 Using \reff{surssol}, we have $\bar w_m(\hat t_ i ,\hat  z_i,\hat d_i ) - \Hc \bar w_m(\hat t_ i ,\hat  z_i,\hat d_i ) \geq
{\delta \over m}$. Therefore, we get
\beqs
\bar \Delta _i & \leq & \underline w (\hat t_i, \hat z_i, \hat d_i)-\bar w_m (\hat t'_i, \hat z_i', \hat d_i')~~\leq ~~ \Hc\underline w (\hat t_i, \hat z_i, \hat d_i)-\Hc\bar w_m (\hat t'_i, \hat z_i', \hat d_i')-{\delta\over m}\;.
\enqs
Sending $i$ to $+\infty$ we get
\beqs
\bar \Delta & \leq & \limsup_{i\rightarrow+\infty}\Hc\underline w (\hat t_i, \hat z_i, \hat d_i)-\liminf_{i\rightarrow+\infty}\Hc\bar w_m (\hat t'_i, \hat z_i', \hat d_i')-{\delta\over m}\\
 & \leq & \Hc\underline w (t_0,  z_0, d_0)-\Hc\bar w_m (t_0, z_0,  d_0)-{\delta\over m}\;,
\enqs
where we used the upper semicontinuity of $\Hc\underline w$ and the lower semicontinuity of   $\Hc\bar w_m$. Since $\underline w$ is upper semicontinuous there exists $a_0\in[0,r_0]$ (with $z_0=(x_0,r_0,p_0,q_0)$) such that $\Hc \underline w(t_0,z_0,d_0)=\underline w(t_0,\Gamma^c(z_0,a_0),d_0)$. Therefore we get the following contradiction
\beqs
\bar \Delta & \leq & \underline w(t_0,\Gamma^c(z_0,a_0),d_0) - \bar w_m(t_0,\Gamma^c(z_0,a_0),d_0) - {\delta\over m} ~~\leq~~\bar \Delta - {\delta\over m}\;.
\enqs
\vspace{2mm}

\ni \textbf{2.} Up to a subsequence we have
\beqs
-\Lc[\hat z_i,\hat d_i,e_i,f_i,M_i] & \leq & 0~~\mbox{  for all }~i\geq 1.
\enqs
Using \reff{surssol} we get
\beq\nonumber
-(e_i-e_i')-  \mu \big(\hat p_i [f_i]_3-\hat p'_i [f'_i]_3\big)-\rho \big(\hat q_i[f_i]_4-\hat q'_i[f'_i]_4\big) & & \\ \nonumber-\eta\big( \hat r_i(\lambda-\hat r_i)[f_i]_2 - \hat r'_i(\lambda-\hat r'_i)[f'_i]_2\big) & & \\ \nonumber
-{1\over 2} \Big(\sigma^2\big(\hat p_i^2[M_i]_{3,3}-\hat p_i'^2[M'_i]_{3,3}\big)+\varsigma^2\big(\hat q_i^2[M_i]_{4,4}-\hat {q_i'}^2[M'_i]_{4,4}\big)& & \\
+2\sigma\varsigma \big(\hat p_i\hat q_i [M_i]_{3,4}-\hat p'_i\hat q'_i [M'_i]_{3,4}\big)+\gamma^2\big(\hat r_i^2 [M_i]_{2,2}-{\hat {r'}_i}^2 [M'_i]_{2,2}\big)\Big)& \leq & -{\delta \over m} \;.\label{ineqishii}
\enq
From \reff{ishii-cond-1}-\reff{ishii-cond-2}, we have
\beqs
e_i ~ = ~ 2(\hat t_i-t_0) & & f_i ~ = ~ 4(\hat z_i-z_0)|\hat z_i-z_0|^2+i(\hat z_i-z_0) \\
e'_i ~ = ~ 2(\hat t'_i-t_0) & & f'_i ~ = ~ 4(\hat z'_i-z_0)|\hat z'_i-z_0|^2+i(\hat z'_i-z_0)
\enqs
and we obtain from \reff{conv Theta i 0} that
\beq\nonumber
-(e_i-e_i')-  \mu \big(\hat p_i [f_i]_3-\hat p'_i [f'_i]_3\big)-\rho \big(\hat q_i[f_i]_4-\hat q'_i[f'_i]_4\big) & & \\-\eta\big( \hat r_i(\lambda-\hat r_i)[f_i]_2 - \hat r'_i(\lambda-\hat r'_i)[f'_i]_2\big) & \xrightarrow[i\rightarrow+\infty]{} & 0\;.\label{first conv}
\enq
Moreover, by \reff{conv Theta i 0} and \reff{ishii-cond-3} , we have using classical arguments
\beqs
\limsup_{i\rightarrow+\infty} \Big(\sigma^2\big(\hat p_i^2[M_i]_{3,3}-\hat p_i'^2[M'_i]_{3,3}\big)+\varsigma^2\big(\hat q_i^2[M_i]_{4,4}-\hat {q_i'}^2[M'_i]_{4,4}\big)\qquad\quad& & \\
 +2\sigma\varsigma \big(\hat p_i\hat q_i [M_i]_{3,4}-\hat p'_i\hat q'_i [M'_i]_{3,4}\big)+\gamma^2\big(\hat r_i^2 [M_i]_{2,2}-{\hat {r'}_i}^2 [M'_i]_{2,2}\big)\Big)& \leq & 0\;.
\enqs
From this last inequality and \reff{first conv} and by sending $i$ to $+\infty$ in \reff{ineqishii} we get $0\leq -{\delta\over m}$, which is the required contradiction.
\vspace{2mm}

\ni$\bullet$ Case 2: we have
\beqs
(\hat t_i,\hat z_i,\hat d_i) & \in & \Dc_k^1\quad \mbox{ for all }~i\geq 1\;.
\enqs
Then we are in the same situation as in the second possibility of Case 1 and we get a contradiction. \ep

\end{document}